\documentclass[journal]{IEEEtran}

\usepackage{float}
\usepackage[caption=false]{subfig}
\usepackage{graphicx}
\usepackage{amsmath}
\usepackage{cite}
\usepackage[dvipsnames]{xcolor}
\usepackage{booktabs}
\usepackage[ruled,vlined]{algorithm2e}
\usepackage{amssymb}
\usepackage[colorlinks = true]{hyperref}
\usepackage{nomencl}

\newtheorem{theorem}{Theorem}
\newtheorem{lemma}[theorem]{Lemma}
\newtheorem{proposition}[theorem]{Proposition}
\newtheorem{definition}[theorem]{Definition}
\DeclareMathOperator{\Tr}{Tr}
\DeclareMathOperator{\eig}{eig}
\DeclareMathOperator{\diag}{diag}
\makenomenclature

\ifCLASSINFOpdf

\else
 
\fi

\begin{document}



\title{Explainable Graph Theory-Based Identification of Meter-Transformer Mapping}

\author{Bilal~Saleem,~\IEEEmembership{Student Member,~IEEE,}
        Yang~Weng,~\IEEEmembership{Senior Member,~IEEE}%
\vspace{-5mm}} 

\maketitle

\begin{abstract}
Distributed energy resources 
are better for the environment but may cause transformer overload in distribution grids, calling for 
recovering meter-transformer mapping to provide situational awareness, i.e., the transformer loading.
\textcolor{black}{The challenge lies in recovering meter-transformer (M.T.) mapping for two common scenarios, e.g., large distances between a meter and its parent transformer or} \textcolor{black}{high} similarity of a meter's consumption pattern to a non-parent transformer's meters. 
Past methods \textcolor{black}{either} 
assume a variety of data as in the transmission grid or ignore \textcolor{black}{the two common scenarios mentioned above.}
Therefore, we propose to utilize the \textcolor{black}{above} observation \textcolor{black}{via} spectral embedding by using the property that 
\textcolor{black}{inter-transformer meter consumptions are not the same and that the noise in data is limited so that all the $k$ smallest eigenvalues of the voltage-based Laplacian matrix are smaller than the next smallest eigenvalue of the ideal Laplacian matrix.}
We also provide a guarantee based on this understanding.
\textcolor{black}{Furthermore, we partially relax the assumption by utilizing location information to aid voltage information for areas geographically far away but with similar voltages.} Numerical simulations on the IEEE test systems and real feeders from our partner utility show that the proposed method correctly identifies M.T. mapping.
\end{abstract}


\IEEEpeerreviewmaketitle


\textcolor{black}{
\printnomenclature
}

\section{Introduction}

With the advent of distributed energy resources (DER), including electric vehicles, the flow of forward and reverse power through the residential electric distribution has increased considerably due to increased consumption and distributed generation. 
\textcolor{black}{Transformers overloaded for a prolonged duration} translate into degraded insulation and a reduced lifespan. Moreover, utilities do not have real-time loading of the transformers to identify \textcolor{black}{overloaded ones}.
Information on real-time transformer loading is \textcolor{black}{essential} but mostly unavailable to utilities as \textcolor{black}{metering} \textcolor{black}{them} increases \textcolor{black}{costs}. Nevertheless, transformer load can also be obtained by summing the loads on \textcolor{black}{downstream smart} meters, where knowledge of meter-transformer \textcolor{black}{(M.T.)} mapping is a prerequisite. 

Unfortunately, M.T. mapping is either unavailable or not updated and obsolete, so it cannot be used for obtaining transformer load.
\textcolor{black}{M.T. mapping information was not required before the advent of DERs~\cite{weng2018big}. Hence it was mostly ignored.}
\textcolor{black}{However, now} electric utilities are \textcolor{black}{intensely} interested in real-time  M.T. mapping. For example, 
\cite{zhang2019meter} identified M.T. mapping using specially designed hardware 
that is very costly and \textcolor{black}{unfeasible} everywhere. Since the literature on the meter-transformer mapping problem is still immature, one idea is to utilize topology identification research to recover \textcolor{black}{such mapping}. 

Topology identification (T.I.) based on data-driven \cite{weng2016distributed,weng2015convexification,yu2017data,weng2014data} graph {structure learning algorithms does not need specially designed hardware. 
For example, \cite{yang2015graphical} uses the Chow-Liu algorithm to recover the entire system topology for radial systems. 
}%
\cite{lauritzen1996graphical} presented an algorithm that works only for decomposable graphs, e.g., a radial system. 
To generalize, researchers developed methods indifferent to topology shape, e.g., radial, looped, or meshed grid. However, such methods require specific voltage probability distributions. For example, \cite{yuan2007model} and \cite{yang2015graphical}
require voltage\textcolor{black}{s} to be multivariate Gaussian and Ising distributed, respectively. 
However, such probabilities are a small subset of the real-world density functions that vary with 
place and time depending on people's usage habits. 
Moreover, such methods need voltage measurements at all system nodes, which is an idealistic assumption for the distribution system domain as service transformers and poles usually do not have any measurements. Moreover, \textcolor{black}{\cite{deka2016estimating} and \cite{cavraro2015data}
require} \textcolor{black}{phasor measurement unit} (PMU) measurements, which are currently not widely available in distribution grids. \textcolor{black}{Similarly, \cite{cavraro2018graph} is based on grid probing using smart inverters, which are not common.}

Some works in T.I. are dedicated to the distribution grid \cite{HaoranLi2021distribution,zhangweng2020topology,liaoWeng2019unbalanced,liaoWeng2018urban,liaoWeng2015distribution}. \textcolor{black}{However}, they require the locations of all switches or the most likely topology, which may be unavailable due to the vast spread of distribution lines~\cite{rudin2011machine}.  
\textcolor{black}{Other} works even require impedance \cite{deka2017structure}, which may be unavailable in the secondary distribution grids.


Advanced metering infrastructure (AMI) data-based {distribution system \textcolor{black}{T.I.} research does not assume voltage at every system node.} For example, \textcolor{black}{\cite{watson2016use} and \cite{luan2015smart}} use smart meter data to determine the topology \textcolor{black}{via estimating the point of connection voltage.} 
Both \cite{luan2015smart} and \cite{watson2016use} assume a fairly accurate prior knowledge of meter-transformer connection information, which is hard to obtain in reality \cite{saleem2020association}. \cite{peppanen2016distribution} uses smart meter voltage, active power, and reactive power measurements via linear regression to determine the topology. However, most residential smart meters do not measure reactive power or power factor. 
\textcolor{black}{\cite{bolognani2013identification} estimates system topology using voltage magnitudes. However, it assumes all lines have the same per unit length inductance to resistance ratios, which is not true. Furthermore,} such methods do not use location information, which is available to most utilities for underground distribution. Even for overhead distribution, utilities can \textcolor{black}{easily} obtain locations \textcolor{black}{using geocoding and} Google Street-view without any field visits. 

Generally, \textcolor{black}{T.I.} is an NP-hard \textcolor{black}{(nondeterministic polynomial time)} problem~\cite{yu2017patopa} and requires \textcolor{black}{a variety} of data~\cite{yu2018patopaem}. Therefore, \textcolor{black}{one can also use clustering to} identify groups of meters belonging to each transformer. \textcolor{black}{Past clustering methods, e.g., $k-$means, `density-based spatial clustering of applications with noise' (DBSCAN), and `balanced iterative reducing and clustering using hierarchies' (BIRCH), fail since they cannot consider two common and \textcolor{black}{challenging} scenarios where a meter’s voltage may be similar to the non-parent cluster, e.g., due to the distance of a meter from its parent transformer or the similarity of a meter's consumption pattern to the non-parent transformer's cluster.
\textcolor{black}{However, going from meter to meter in voltage feature space leads to the correct meter cluster. It is because, within a transformer secondary feeder, a meter’s voltage is usually more similar to a neighboring meter’s voltage than the mean voltage of all meters in the feeder.}}

\textcolor{black}{In this paper}, we propose to utilize the \textcolor{black}{above} observation by employing spectral embedding with an explainable working mechanism.
Although we can use the spectral embedding-based method, the next challenge is to \textcolor{black}{provide a theoretical guarantee for the method.}
Therefore, we propose a proof of performance guarantee 
under the following reasonable assumption.
\textcolor{black}{The consumptions of inter-transformer meters are not the same, and the noise in voltage data is limited so that all the $k$ smallest eigenvalues of the voltage-based Laplacian matrix are smaller than the next smallest eigenvalue of the ideal Laplacian matrix.}

However, such an assumption \textcolor{black}{may} restrict the algorithm from working \textcolor{black}{in} areas where 
meter voltages are \textcolor{black}{alike} due to similar consumptions or \textcolor{black}{very high penetration of} DERs,
so the next challenge is to relax the assumption. Therefore, 
we partially relax the assumption on the Laplacian matrix spectrum by utilizing location information in addition to the voltage information. 
\textcolor{black}{Such a method also works for looped or meshed feeders, and it does not need a specific voltage probability density \cite{von2007tutorial}. Also, it does not need an initial version of meter-transformer} \textcolor{black}{mapping information as an input.} 



Numerical experiments are carried out 
on the standard distribution testbeds, e.g., IEEE $123$-bus and our partner utility's local grid with \textcolor{black}{a} half-million customers. The \textcolor{black}{results validate the assumption and demonstrate} that proposed solutions accurately segment the smart meter data to identify \textcolor{black}{meter-transformer mapping.}


The rest of this paper is organized as follows: Section II shows problem modeling. Section III presents algorithms. Section IV provides guarantees. 
Section V includes location information. \textcolor{black}{Finally, s}ection VI evaluates performance, and Section VII concludes the paper.

\begin{figure}[t]
\vspace{-1mm}
\centering
\includegraphics[width = 0.3 \textwidth]{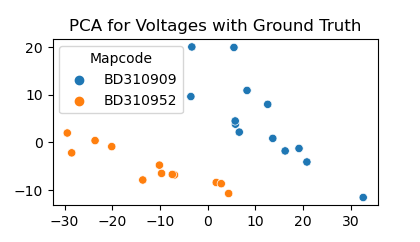}
\vspace{-3mm}
\caption{shows a case of meter-transformer connectivity from a local utility in Arizona, showing 
the ground truth and voltage ``nearness'' of meters \textcolor{black}{to} non-parent clusters.}
\label{fig_meter_near_wrong_cluster}
\vspace{-5mm}
\end{figure}

\vspace{-2mm}

\section{Problem Modelling}


\textcolor{black}{To define the meter-transformer (M.T.) mapping identification method, we describe time-series voltage data given by smart meters.} For example, the latitude-longitude pairs in radians for $N$ smart meters $\mathbf{l}^1,\cdots,\mathbf{l}^N \in \mathbb{R}^{2 \times 1}$ are stored as row vectors in matrix $L \in \mathbb{R}^{N \times 2}$, where $\mathbb{R}$ represents the set of real numbers. The voltage time-series  with $T$ timeslots for $N$ smart meters $\mathbf{v}^1,\cdots,\mathbf{v}^N \in \mathbb{R}^{T \times 1}$ are stored as row vectors in matrix $V \in \mathbb{R}^{N \times T}$. In addition to smart meters, we assume that there are $k$ transformers forming $k$ clusters of smart meters in the distribution grid. 
\textcolor{black}{Also, transformer locations are available.}
\textcolor{black}{$C_j$} represents a set of indices of all smart meters in the $j^{th}$ cluster. 
A smart meter $i \in \{1, \cdots, N\}$ is uniquely present in a cluster $j \in \{1, \cdots, k\}$ that is supplied by a common transformer. There exists a many-to-one mapping $f:i \rightarrow  j$. 

For correlating these variables, a distribution system is characterized by buses $\mathcal{V} = {1, 2, \cdots, N}$ and by branches $\mathcal{E} = {(i,i'),i,i' \in \mathcal{V}}$. The voltage measurement data at bus $i$ and time $t$ \textcolor{black}{is} represented as the magnitude of the instantaneous voltage at bus $i$ in per-unit  
$|v_i(t)| \in \mathbb{R}$.
The \textcolor{black}{meter voltage} measurements in $\mathbf{v}^i$ are steady-state voltages over a period according to utility collection speed. 
We define the problem below.
\begin{itemize}
\item Problem: identify smart meter to transformer connectivity 
\item Given: smart meter voltage data $V$ and the smart meter location data $L$.
\item Find the \textcolor{black}{M.T.} mapping $f:i \rightarrow  j$.
\end{itemize}
\textcolor{black}{We do not consider active
power consumption since they are dependent on the consumers and are not significantly affected by the \textcolor{black}{M.T.} mapping. $f$ is the desired M.T. mapping.}

\vspace{-2mm}

\section{\textcolor{black}{Clustering-Based M.T. Mapping Identification}}

\label{sec:DBSCAN_Kmeans_BIRCH_SpecClust}

\textcolor{black}{As discussed in the Introduction, meter-transformer (M.T.) mapping is needed to identify overloaded transformers. \textcolor{black}{However, we cannot use voltage measurements from transformers to identify the M.T. mapping since they are \textcolor{black}{unavailable} at many utilities. Nonetheless,} we obtain transformer-based meter clusters and find the nearest transformer to each cluster. Such an approach is better than obtaining the nearest transformer to each meter.} \textcolor{black}{It is because} \textcolor{black}{real cluster-centers are more likely to be closer to their parent transformer than individual smart meters, which might be at the cluster’s \textcolor{black}{boundary.} However, the challenge lies in obtaining the meter clusters supplied by a common transformer.} 
Below is an evaluation of clustering methods for our purpose.



\vspace{-2mm}

\subsection{Metric Evaluation for Clustering Algorithm Design}

For clustering data, three categories are popular in data mining.  
One is to consider the combined properties of clusters, e.g., minimize the sum of distances within each cluster ($k-$means). The second category is to \textcolor{black}{set} bounds on clusters, e.g., maximum diameter for clusters (BIRCH). \textcolor{black}{Finally, t}he third category is based on the orthogonal eigenvectors to separate the smart meter clusters (spectral clustering). 


\textcolor{black}{We analyze the representative methods from the three classes by analyzing their suitability for power distribution data, including voltage and location data.}

\subsubsection{\texorpdfstring{$k-$}{k-}means}

One idea for clustering is to consider the within-cluster properties of all the members in a group. For example, $k-$means forms $k$ groups of smart meters.
It aims at minimizing the squared error loss. The algorithm works by alternating centroid computation and cluster adjustment steps.

\textit{Drawbacks for Power Data}: 
$k-$means tends to produce spherical clusters due to the minimization of within-cluster distances since a sphere has the minimum within-cluster distances for the same area compared to any other shape. A real transformer secondary circuit may have an irregular shape due to irregular street shapes. 

\subsubsection{BIRCH for Maximum Cluster Distance}
Instead of considering within-cluster parameters as in $k-$means, one can also \textcolor{black}{set a bound} on the extreme points within a cluster, where the BIRCH algorithm is well-known. \textcolor{black}{Each} point is assigned to the nearest-\textcolor{black}{centered} subcluster. 

\textit{Drawbacks for Power Data}: The \textcolor{black}{secondary distribution from a transformer} may have an irregular shape depending on the shape of the street. Setting a hard limit on the radius requires all transformer secondaries to be limited to the sphere. For example, both long and short feeders are popular in the power domain. Therefore, it is unwise to have a hard limit on cluster diameter. 

\subsubsection{Proposed Spectral Theory-Based Approach}
\begin{figure}[t]
    \centering
    \includegraphics[width = 0.48 \textwidth]{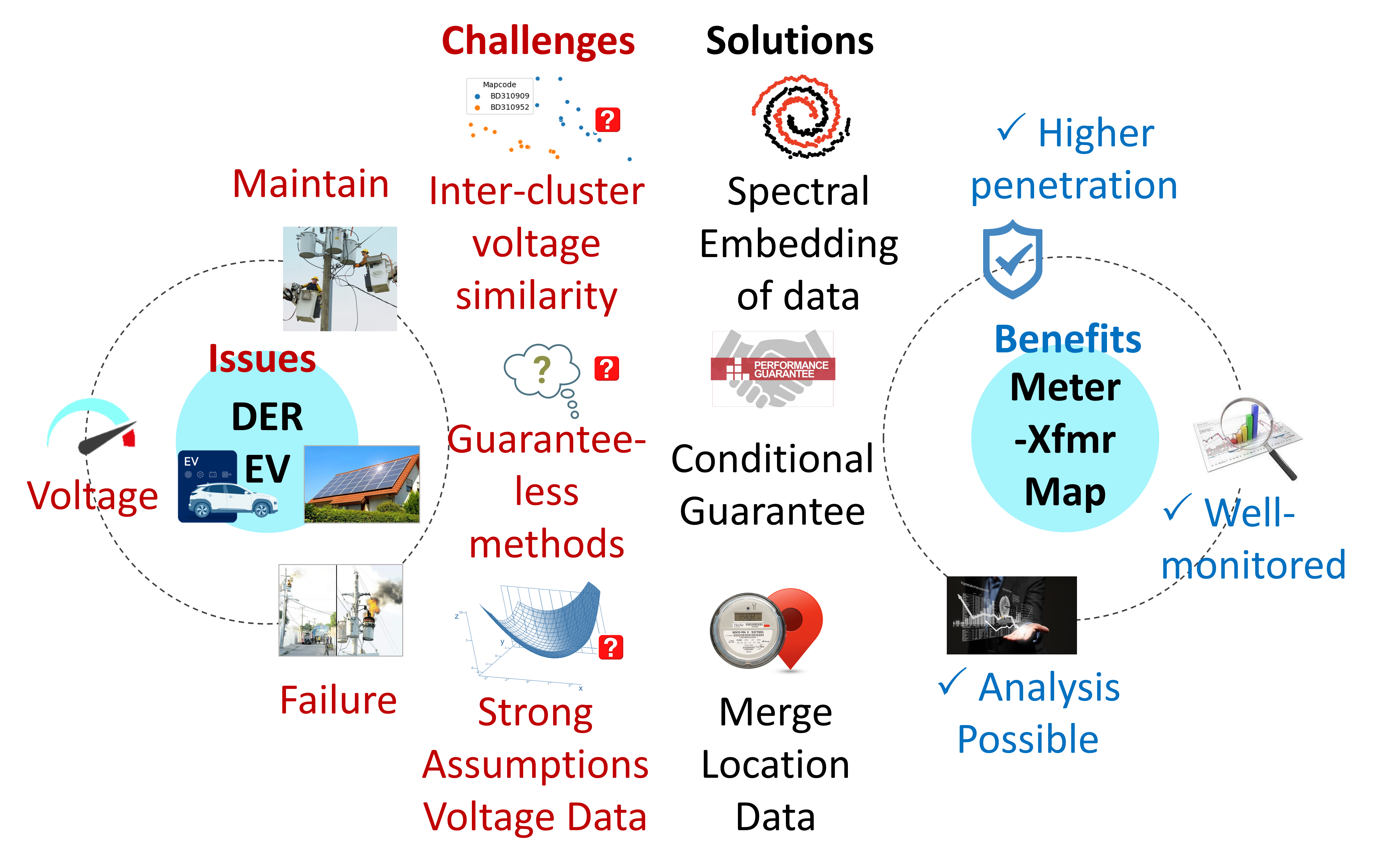}
    \vspace{-2mm}
    \caption{summarises the \textcolor{black}{topology-related challenges and the contributions of the proposed methods to resolve the challenges.}}
    \label{fig:Overview_diagram}
    \vspace{-2mm}
\end{figure}

The two approaches discussed above are based on within-cluster distance and maximum cluster diameter. Another idea is to use eigenvectors of the graph Laplacian matrix to separate the clusters. Such an approach is widely popular in graph theory for graph separation applications.\\
\textit{Suitability for Power Data:} The graph Laplacian matrix \textcolor{black}{is} estimated from the voltage data. The spectral theory-based approach does not need clusters of a definite shape or size.
Moreover, \textcolor{black}{ground truth} clusters can be identified via a provable guarantee. Fig.~\ref{fig:Overview_diagram} shows an overview of the challenges, the proposed spectral clustering-based solutions, and the obtained benefits. \textcolor{black}{Therefore, the proposed approach is better than other clustering approaches for M.T. mapping identification due to its suitability for power system data.}

\vspace{-2mm}

\subsection{Spectral Clustering}

\textcolor{black}{In the remaining part of this section and Section~\ref{Sec:Spectral_Cluster_Guarantee}, we demonstrate spectral clustering working and guarantee solely using voltage data under the assumption mentioned in the Introduction. Finally, in Section~\ref{Sec:Co-regularized_multi-view_Spectral_Clustering}, we will demonstrate the use of location data to aid voltage data for partially relaxing the assumption. Fig.~\ref{fig_meter_near_wrong_cluster} shows a challenging scenario for recovering meter-transformer mapping where some meter voltages are more similar to the center of the non-parent transformer cluster. Below, we show the motivation of the similarity matrix, which is needed for spectral clustering.}

\begin{figure}[ht]
    \centering
    \includegraphics[width = 0.35 \textwidth]{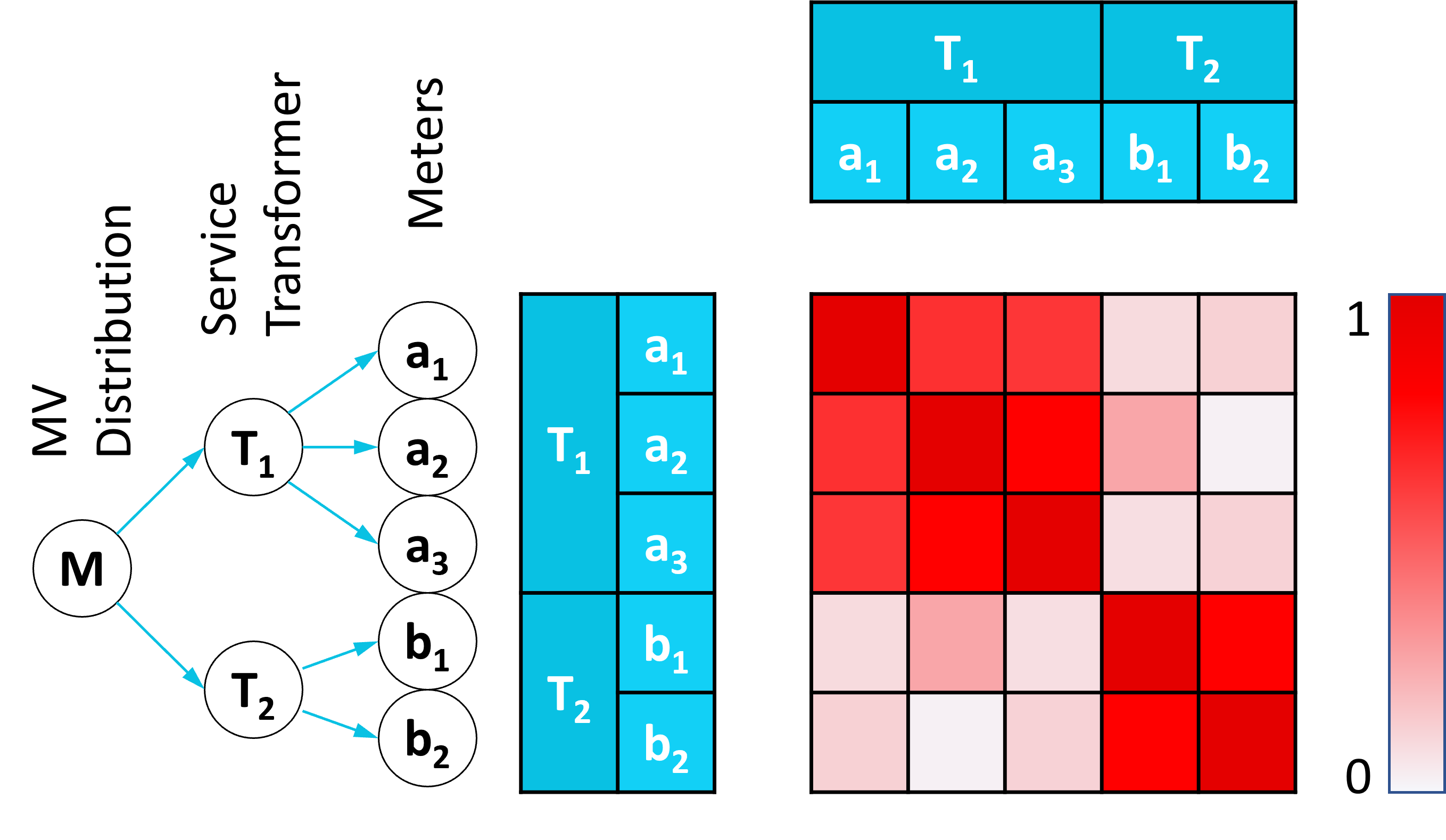}
    \vspace{-3mm}
    \caption{shows a similarity matrix with submatrices formed by \textcolor{black}{two} transformers' secondaries.}
    \label{fig:spectral_Clustering_Introduction}
\end{figure}

\textcolor{black}{A similarity matrix is a square matrix $M \in  \mathbb{R}^{N \times N}$,} where the rows and columns correspond to the smart meters. Each entry $m_{ii'}$ represents the similarity metric between the smart meters $i$ and $i'$, as shown in Fig.~\ref{fig:spectral_Clustering_Introduction}. The similarity metric is constructed as $m_{ii'} = \exp \left( \frac{-d^2(\mathbf{v}^i,\mathbf{v}^{i'})}{\sigma^2}\right)$ for $i \neq i'$ and $m_{ii} = 1$, 
\textcolor{black}{where $d(\mathbf{v}^i,\mathbf{v}^{i'})$ is the Euclidean distance.}
$\sigma$ is a scale parameter. It can be seen that $\sigma$ enhances or \textcolor{black}{diminishes} the effect of $d(\mathbf{v}^i,\mathbf{v}^{i'})$. 

\textcolor{black}{By forming a similarity matrix,} we reduced the dimensions of the input data from $N \times T$, where $T$ is the number of timestamps of the voltage data and $T > N$, to $N \times N$. In doing so, we also neglected the unwanted information, such as the mean value of smart meter voltage. 
Since a power grid \textcolor{black}{is} modeled as a graph, the similarity matrix $M$ is a scaled approximation of the weighted graph adjacency matrix \textcolor{black}{$E \in \mathbb{R}^{N \times N}$.}
$E_{ii'}$ is the weight of the edge between nodes $i$ and $i'$ in graph $G$. 
\textcolor{black}{An unweighted graph adjacency matrix is a binary matrix. $E_{ii'}$ is one if there is an edge $\{i,i'\}$ otherwise, zero.} 

Given \textcolor{black}{the weighted} graph adjacency matrix or similarity matrix, we use our method to cluster the vertices of an unknown graph into $k$ clusters, where $k$ is the number of desired clusters. For distribution grid graphs, the vertices are smart meters, and $k$ is the number of transformers in an area. 
Also, higher edge weights indicate a higher degree of similarity of smart meter voltages.









\label{sec:Spectral_Clustering_for_Multiple_Clusters}

\textcolor{black}{To obtain $k$ clusters, one can minimize the sum of edge weights between $k$ clusters
$cut(C_1,C_2,\cdots,C_k) := \sum_{i=1}^{k} cut(C_i,\overline{C_i})$, where $cut(C_i,\overline{C_i})$ is defined as the sum of inter-cluster similarity matrix elements \textcolor{black}{$m_{ij}$} corresponding to cluster $C_i$ data points.}
However, such minimization will result in all meters into a single cluster, \textcolor{black}{as expected}. Therefore, we divide $cuts$ by the number of vertices in each cluster to obtain $\mathcal{Q} (C_1,C_2,\cdots,C_k) = \sum_{i=1}^{k} \frac{cut(C_i,\overline{C_i})}{|C_i|}$.
We can rewrite the above relation in terms of the \textcolor{black}{unnormalized graph Laplacian $\mathcal{L} = D-M$, where \textcolor{black}{$M$} is the similarity matrix and \textcolor{black}{$D$} is the diagonal degree matrix with $d_{ii} = \sum_{j} m_{ij} $. For example, i}f we define an indicator vector $\mathbf{h_i}$ as follows 
\begin{align}
 h_{i,j} =
    \begin{cases}
        1/\sqrt{|C_j|} & \text{if $i \in C_j$,}\\
        0 & \text{if $i \notin\overline{C_j}$,}
    \end{cases}
    \label{Eq:h_definition}
\end{align} 
then $\mathbf{h_i}^T \mathcal{L} \mathbf{h_i} = \left.{cut(C_i, \overline{C_i})} \middle/ {|C_i|} \right.$ and $\mathcal{Q} (C_1,C_2,\cdots,C_k) = \sum_{i=1}^{k} \mathbf{h_i}^T \mathcal{L} \mathbf{h_i}$.
Let $H \in \mathbb{R}^{n \times k}$ be such that the $k$ indicator vectors $\mathbf{h_i}$ are the columns of $H$. We see that $\mathbf{h_i}^T \mathcal{L} \mathbf{h_i} = \left ( H^T \mathcal{L} H \right )_{ii}$. Moreover, $\mathcal{Q} (C_1,C_2,\cdots,C_k) = \sum_{i=1}^{k} \left ( H^T \mathcal{L} H \right )_{ii} = \Tr \left ( H^T \mathcal{L} H \right ),$ where $\Tr$ is the trace of a matrix. So the combined optimization problem is defined as below. 
\begin{align*}
&\min_{C_1,C_2, \cdots,C_k} \Tr \left ( H^T \mathcal{L} H \right ) \\
& \text{subject to $H^T H = I$, where $H$ is defined from Eq.~\ref{Eq:h_definition}}. 
\end{align*} \textcolor{black}{The constraint $H^T H = I$ specifies that each data point belongs to a single cluster. Minimizing $\Tr \left ( H^T \mathcal{L} H \right )$ is a discrete optimization problem as $h_{i,j}$ takes only two values, \textcolor{black}{and} it is an NP-hard problem. Therefore,} we relax the problem by allowing $H$ to take on continuous values instead of the discrete values.
\begin{align}
&\min_{H \in \mathbb{R}^{N \times k} } \Tr \left ( H^T \mathcal{L} H \right ) \label{Eq:Regular_Spectral_Clustering} \\
& \text{subject to $H^T H = I$}. \nonumber
\end{align}
\begin{theorem}
\label{theorem_second_eigenvalue}
The solution to the above constraint minimization problem \textcolor{black}{is the $H$ comprising $k$ eigenvectors corresponding to the $k$ smallest eigenvalues of $\mathcal{L}$.}
\end{theorem}
\begin{IEEEproof}
\textcolor{black}{See Appendix~\ref{Appendix_Proof_of_Theorem_1} for proof.}
\end{IEEEproof}
\textcolor{black}{Therefore, the problem becomes an eigendecomposition problem. Moreover, such a solution satisfies} $H^T H = I$, since normalized eigenvectors are orthogonal to each other with unity magnitude.\\
Finally, we use the $k-$means++ algorithm to convert the continuous values of $H$ into discrete clusters. \textcolor{black}{The spectral embedding makes the irregular-shaped clusters into regular shapes, so $k-$means++ can identify the correct clusters, as we will elaborate in Section~\ref{Sec:Spectral_Cluster_Guarantee}.} The steps for the algorithm are shown below.\\
Given $N$ smart meters voltage time series $\mathbf{v}^1,\cdots,\mathbf{v}^N$. We cluster them into $k$ transformer secondary clusters as follows:

\begin{enumerate}
\item Form the affinity matrix $M \in \mathbb{R}^{n \times n}$ as defined earlier.
\item Define $D$ to be a diagonal matrix with $d_{ii} = \sum_{j=1}^{n} m_{ij} $ and construct the graph Laplacian matrix $\mathcal{L} = D-M$. The off-diagonal elements of the Graph Laplacian represent the similarity between the two nodes \textcolor{black}{$m_{ii'}$}. In contrast, the diagonal elements represent the total similarity of a node \textcolor{black}{$d_{ii}$.}
\item Select the desired number of groups $k$ as the number of transformers.
\item Find $\mathbf{x}_{(1)},\cdots,\mathbf{x}_{(k)},$ the eigenvectors corresponding to $k$ smallest eigenvalues of $\mathcal{L},$ and form the matrix $X = \left[ \mathbf{x}_{(1)},\cdots,\mathbf{x}_{(k)} \right] \in \mathbb{R}^{n \times k}$. Since eigenvectors are orthogonal to each other, doing so will further distance the points belonging to different clusters.
\item Treat each row of $X$ as a point in $\mathbb{R}^k$ and cluster via $k-$means$++$. Due to orthogonalizing using the eigenvectors, the data points belonging to separate clusters are almost orthogonal to each other, i.e., they have approximately right angles at the origin with respect to (w.r.t.) each other\textcolor{black}{, so that $k-$means++ can cluster well.}
\item For the matrix $X$'s rows assigned to cluster $C$, the original corresponding points $s_i$ are present in cluster $C$ \cite{ng2002spectral}.
\end{enumerate}

So far, we have shown a method to cluster meters into groups based on service transformer, although some meters have similar voltage to the non-parent group. The question that arises \textcolor{black}{concerns the} guarantee for \textcolor{black}{it}. \textcolor{black}{The next section provides} a guarantee for our proposed method and the \textcolor{black}{rational} assumption needed for the guarantee.

\section{\textcolor{black}{Theoretical Guarantee for Spectral Clustering Robustness}}
\label{Sec:Spectral_Cluster_Guarantee}

\textcolor{black}{In this section, we provide a guarantee for spectral clustering under \textcolor{black}{the assumption that} 
\textcolor{black}{meter consumptions are not the same, and the noise in data is limited so that all the $k$ smallest eigenvalues of the voltage-based Laplacian matrix are smaller than the $(k+1)$-th smallest eigenvalue of the ideal Laplacian matrix.}
In our analysis below, we go from simple to more realistic scenarios.}
\vspace{-5mm}
\subsection{\textcolor{black}{Simple Scenario (Laplacian as a Block-Diagonal Matrix):}}
\textcolor{black}{This subsection simplifies the scenario by a stronger assumption that voltages from smart meters of different transformers are independent. However, such an assumption is strong and not very useful. Therefore, in the next subsection, we will replace it with a weaker and more useful assumption, i.e., the $k$ smallest eigenvalues of the voltage-based Laplacian matrix are smaller than the $(k+1)$-th smallest eigenvalue of the ideal Laplacian matrix.}
\textcolor{black}{Following the strong assumption on voltage-independency, as described above, we have that} the similarity matrix consists of values $0 \leq m_{ij} \leq 1$. Moreover, the elements
$m_{ij} > 0$ indicate meter $i$ and meter $j$ belong to the same transformer; otherwise, \textcolor{black}{different transformers.}
Let $n_i$ be the number of smart meters in the $i-$th cluster. 

Furthermore, \textcolor{black}{upon permutation, meters can be assumed to be in consecutive columns/rows of the similarity matrix without loss of generality}. Therefore, the similarity matrix $M$ has a block-diagonal structure. Mathematically, $M = \diag(M_{n_1}, M_{n_2}, \cdots, M_{n_k})$,  \textcolor{black}{where $\diag(\cdots)$ represents a diagonal matrix or block-diagonal matrix with elements $(\cdots)$.} Each diagonal block is a square submatrix $M_{n_i}$ with dimensions $n_i$. \textcolor{black}{The Laplacian Matrix $\mathcal{L} = D-M$, \textcolor{black}{where $D$ is the diagonal degree matrix,} will also be block-diagonal.}
\textcolor{black}{\subsubsection{Ideal case: Zero impedance of wire from meter-transformer}
In such a case, the meters downstream of a transformer will have the same voltage. Moreover, its corresponding entry in $M$ will be $m_{ij} = 1$. We also assume the entries in $M$ for different transformers are zero, i.e., $m_{ij} = 0$.}


\begin{lemma}
For an $n-$dimensional \textcolor{black}{square} matrix $\left(nI_n-J_n\right)$ with diagonal values $(n-1)$ and off-diagonal values $(-1)$, the smallest eigenvalue is $0$ with the respective eigenvector as the $n$-dimensional all-ones vector $\mathbf{1}_n$, all other eigenvalues are $n$.
\end{lemma}

\begin{IEEEproof}
The lemma is obvious, so we skip the proof.
\end{IEEEproof} 
\textcolor{black}{Lemma 2 relates to an unweighted graph adjacency matrix but not a similarity matrix. On the other hand, lemma 3 considers a block-diagonal matrix.}
\begin{lemma}
The eigenvalues of a block diagonal matrix are the union of the eigenvalues of constituent diagonal blocks. Also, the eigenvectors of the matrix are the union of the eigenvectors of diagonal blocks padded appropriately with zeros based on the location of the respective diagonal block.
\end{lemma}

\begin{IEEEproof}
Consider $B = diag([ B_{n_i}]_{i=1,\cdots,k})$ and $\mathbf{\hat{x}}$ a zero-padded eigenvector of $B_{n_i}$. Remaining proof is obvious.
\end{IEEEproof}

\textcolor{black}{The Laplacian matrix $\mathcal{L}$ 
is an $N \times N$ matrix corresponding to the $N$ smart meters. \textcolor{black}{The diagonal value of a similarity matrix $M$ is not useful for forming meter clusters for meter-transformer mapping.} Therefore, we construct the Laplacian matrix $\mathcal{L} = D-M$ from the similarity matrix to neglect the effect of the diagonal entries of the similarity matrix. The Laplacian matrix has exactly $k$ eigenvalues that are zeros, where $k$ is the number of smart meter clusters that is the same as the number of transformers. \textcolor{black}{ Let $X := \{\mathbf{x}_1, \mathbf{x}_2, \cdots , \mathbf{x}_k\} \in \mathbb{R}^{n \times k}$ be the matrix containing eigenvectors of $\mathcal{L}$ corresponding to eigenvalue $0$.} The $i-$th eigenvector \textcolor{black}{$\mathbf{x}_i$ can} have ones at indices corresponding to the $n_i$ meters of the $i-$th cluster.}

\textcolor{black}{However, an important point needs to be made clear. As $0$ is a repeated eigenvalue, the eigenvectors can be any \textcolor{black}{$k$} orthogonal vectors covering the same subspace as $k$ \textcolor{black}{zero-padded} eigenvectors \textcolor{black}{from the union of the individual blocks.} \textcolor{black}{In other words, $X$ can be replaced by $XQ$ for any orthogonal matrix $Q \in \mathbb{R}^{k \times k}$ ($Q^TQ=QQ^T=I$)}. \textcolor{black}{Thus, a guarantee cannot be provided for the individual eigenvectors. However, a guarantee \textcolor{black}{may} be provided for the subspace (any linear combination of the $k$ eigenvectors) corresponding to eigenvalue $0$.}}

\begin{theorem}
Let $X := \{\mathbf{x}_1, \mathbf{x}_2, \cdots , \mathbf{x}_k\} \in \mathbb{R}^{n \times k}$ be the matrix containing eigenvectors of $\mathcal{L}$ corresponding to eigenvalue $0$. 
\begin{enumerate}
\item The row vectors  of $X$ corresponding to the data points of the same cluster are equal. 
\item The row vectors of $X$ corresponding to different clusters are orthogonal.
\end{enumerate}
\end{theorem}


\begin{IEEEproof}
\textcolor{black}{See Appendix~\ref{Appendix_Proof_of_Theorem_4} for proof.}
\end{IEEEproof}

\begin{proposition}
$k-$means++ algorithm \textcolor{black}{identifies} true meter clusters using the matrix $X$.
\end{proposition}

\begin{IEEEproof}
The proposition is obvious by considering the initialization procedure of $k-$means++ algorithm. 
\end{IEEEproof}


\vspace{-3mm}
\textcolor{black}{\subsubsection{Non-Ideal Case: Relaxation of the ideal case}
\label{Laplacian_as_a_Block_Diagonal_Matrix}
In such a case, we allow continuous values for the entries $0 \leq m_{ij} \leq 1$.}


\begin{proposition}
The block-diagonal Laplacian matrix has $0$ eigenvalue repeated $k$ times.
\end{proposition}

\begin{IEEEproof}
A proof can be done by lemma 3, and that $0$ is an eigenvalue of each diagonal block $\mathcal{L}_{n_i}= D_{n_i} - M_{n_i}$.
\end{IEEEproof}



\textcolor{black}{According to \textcolor{black}{proposition 6}, like the simple scenario, the matrix of eigenvectors $X := \{\mathbf{x}_1, \mathbf{x}_2, \cdots , \mathbf{x}_k\} \in \mathbb{R}^{n \times k}$ of $\mathcal{L}$ corresponding to eigenvalue $0$ has equal row vectors for data points of the same cluster. Moreover, the row vectors of $X$ corresponding to different clusters are orthogonal. \textcolor{black}{Thus, according to proposition 5}, the $k-$means++ algorithm \textcolor{black}{identifies} the true meter clusters.}










\textcolor{black}{In the case of similar inter-transformer voltages, e.g., high penetration of behind-the-meter photovoltaic (PV), the simple scenario does not remain valid since the assumption of block diagonal matrix is no longer valid. However, we will resolve such a scenario in the next subsection.}

\textcolor{black}{\textbf{Exceptional Scenario} In the case of extremely similar inter-transformer voltages, e.g., having a Pearson correlation coefficient equal to $1$, any voltage-based method can not work. Therefore, we exclude such a scenario from our scope by \textcolor{black}{using the first part of our assumption, i.e., meter consumptions are not the same.} Thus, the voltages can not be exactly similar due to the different drops across transformers and wires.} \textcolor{black}{So far, the Laplacian matrix as a block diagonal matrix is considered. In the next subsection, the general scenario is presented.}

\vspace{-3mm}

\subsection{General Scenario:}


In the general scenario, \textcolor{black}{we replace the strong assumption from the simple scenario, i.e., we no longer assume voltage independence among meters supplied by different transformers. However, we assume the weaker assumption that \textcolor{black}{the $k$ smallest eigenvalues of the voltage-based Laplacian matrix are smaller than the $(k+1)$-th smallest eigenvalue of the ideal Laplacian matrix.} Therefore, the Laplacian matrix is no longer block-diagonal.
Moreover,} there may also be non-zero repeated eigenvalues of the \textcolor{black}{symmetric} Laplacian matrix $\mathcal{L}$. In such a case, the eigenvectors are not unique. In fact, any orthogonal transformation of the eigenvectors will yield a set of eigenvectors for the same eigenvalue. However, the subspace spanned by the eigenvectors of an eigenvalue is unique. 
Therefore, the concept of simple invariant subspaces is useful. \textcolor{black}{Below a brief introduction to the simple invariant subspaces is provided.}

\textcolor{black}{Consider $\lambda_1$ a repeated eigenvalue in $\mathcal{L}$. Let $X_1$ be a matrix whose columns are a set of eigenvectors \textcolor{black}{with unit magnitudes} of the Laplacian matrix $\mathcal{L}$ for the eigenvalue $\lambda_1$. The columns of matrix $X_1$ span a subspace, denoted by $\mathcal{R}(X_1)$, such that $\mathcal{LR}(X_1) \subseteq \mathcal{R}(X_1)$.
Such a subspace $\mathcal{R}(X_1)$ is a simple invariant subspace of the matrix $\mathcal{L}$. For all $\mathbf{x} \in \mathcal{R}(X_1)$, \textcolor{black}{it’s reflected that} $\mathcal{L} \mathbf{x} = \lambda_1  \mathbf{x}$. Moreover, $\mathcal{L} X_1 = \lambda_1 X_1$, and $\lambda_1 I = X_1^{T} \mathcal{L} X_1$ is the corresponding eigenvalue since $X_1$ is \textcolor{black}{an orthogonal} matrix (i.e., $X_1^{T} X_1 = X_1 X_1^{T} = I$).
Therefore, the subspace $\mathcal{R}(X_1)$ is \textcolor{black}{also the eigenspace} of the matrix $\mathcal{L}$ \textcolor{black}{corresponding to the eigenvalue $\lambda$}.}

\textcolor{black}{Let $Y_1$ be a matrix of all eigenvectors \textcolor{black}{with unit magnitudes} of $\mathcal{L}$, except \textcolor{black}{those} corresponding to eigenvalue $\lambda_1$. Since $\mathcal{L}$ is symmetric, the columns of $Y_1$ are orthogonal to the columns of $X_1$. Therefore, $Y_1^{H} \mathcal{L} X_1=0$, which is a necessary condition for $\mathcal{R}(X_1)$ to be an invariant space. By transposing on both sides, we get $X_1^{H} \mathcal{L}^H Y_1= X_1^{H} \mathcal{L} Y_1 = 0$, \textcolor{black}{which indicates} that both $X_1$ and $Y_1$ are simple invariant subspaces of $\mathcal{L}$ so that $\mathcal{L}Y_1 = Y_1 \Lambda_1$
, where $\Lambda_1$ is a diagonal matrix of all eigenvalues of the Laplacian matrix $\mathcal{L}$, \textcolor{black}{except $\lambda_1$}. To \textcolor{black}{summarize}, a simple invariant subspace is a more generalized version of \textcolor{black}{eigenspaces}. \textcolor{black}{For example,} for non-repeated eigenvalues, the corresponding 
eigenspaces are \textcolor{black}{straight, infinitely long} lines along the direction of the respective eigenvectors.}

\textcolor{black}{For the general scenario, we consider the Laplacian matrix $\mathcal{\tilde{L}}$ as a perturbed version of the block-diagonal Laplacian matrix $\mathcal{L}$.} Let $\mathcal{\tilde{L}} = \mathcal{L} + d\mathcal{L}$, where $d\mathcal{L}$ is due to the impact of the non-idealized environment. Moreover, let $(\lambda_1,X_1)$ be the eigenvalue and the corresponding matrix whose columns are the eigenvectors of $\mathcal{L}$, $\mathcal{R}(\tilde{X}_1)$ be an \textcolor{black}{approximation} of $\mathcal{R}(X_1)$, 
\textcolor{black}{and $\tilde{P}_1 = \tilde{X}_1^T \mathcal{L} \tilde{X}_1$} be the approximation of $\lambda I$.

\textcolor{black}{In order to have $\mathcal{R}(\tilde{X}_1)$ close to $\mathcal{R}(X_1)$, \textcolor{black}{the norm of \textcolor{black}{$R := \mathcal{L} \tilde{X}_1 - \tilde{X}_1 \tilde{P}_1$} needs to be small enough.} Moreover, we need a sufficient gap between $\lambda_1$ and the other eigenvalues of $\mathcal{L}$. }
\textcolor{black}{The gap between eigenvalues is necessary to ensure $\mathcal{L} \tilde{X}_1$ is less affected by the other eigenvalues of $\mathcal{L}$.}

\begin{definition}
the set of all eigenvalues of a matrix A is defined as $\eig(A)$.
\end{definition}

\begin{theorem}
Let $\mathcal{L}$ be symmetric. Also, let the columns of $X_1$ form an orthogonal basis for the simple invariant subspace of $\mathcal{L}$ for the eigenvalue $\lambda$. Moreover, let $\tilde{X}_1$ be the approximation of $X_1$, and $\tilde{P}_1 = \tilde{X}_1^T \mathcal{L} \tilde{X}_1$ become the approximation of $\lambda I$. \textcolor{black}{Moreover, let $\eig(\tilde{P}_1) \subseteq [a,b]$.} Let  $R = \mathcal{L}\tilde{X}_1 - \tilde{X}_1 \tilde{P}_1$ be the residual of the approximation $\tilde{X}_1$. Moreover, let 
\textcolor{black}{$sep \left([a,b],\eig(\Lambda_2)\right)>\delta$, where $sep(\cdot,\cdot)$ is the minimum distance over all elements of the two sets.} Then,
$$||\tan {\Theta} [\mathcal{R}(\tilde{X}_1),\mathcal{R}(X_1)] || \leq \frac{||R||}{\delta}.$$
\end{theorem}

\begin{IEEEproof}
See Appendix~\ref{Appendix_Proof_of_Theorem_8} for proof.
\end{IEEEproof}

\textcolor{black}{\textit{Intuition for $\delta$ (need for the assumption):}
The gap between $\lambda$ and the set $\Lambda$ is known as the eigengap. It is a measure of the stability of the invariant subspaces against the perturbation of the Laplacian $\mathcal{L}$ matrix. For example, if the perturbation is very large, or the eigengap is too small, then $\eig(\tilde{P}_1)$ might overlap $\eig(\Lambda)$, and the guarantee does not exist. \textcolor{black}{As we will demonstrate below, t}he eigengap is larger than the gap between the eigenvalues of the perturbed space $\eig(\tilde{P}_1)$. If only the perturbation is less than the eigengap\textcolor{black}{,} we have the guarantee that the difference between the two subspaces is bounded. } 

\begin{theorem}
The gap between $\lambda$ and the set $\Lambda$ (eigengap) is larger than the gap between the eigenvalues of the perturbed space $\eig(\tilde{P}_1)$.
\end{theorem}

\begin{IEEEproof}
See Appendix~\ref{Appendix_Proof_of_Theorem_9} for proof.
\end{IEEEproof}









As discussed in the Introduction, meter-transformer mapping is a prerequisite to \textcolor{black}{resolve} the challenges in the distribution grid due to distributed energy resources (DERs), including electric vehicles. \textcolor{black}{So far, the spectral clustering algorithm and its guarantee with the required assumption to resolve the challenge have been discussed.} 
However, the assumption \textcolor{black}{may} limit the applicability of the algorithm. \textcolor{black}{Therefore, in the next section, partial relaxation of the condition on the Laplacian matrix will be addressed by utilizing meter location information widely available to utilities.}
\vspace{-3mm}
\section{Co-regularized Multi-view Spectral Clustering }
\label{Sec:Co-regularized_multi-view_Spectral_Clustering}


\textcolor{black}{To generalize the \textcolor{black}{algorithm's applicability}, we focus on relaxing the condition. \textcolor{black}{In this section, we use meter location information with the voltage information} to partially relax \textcolor{black}{the assumption that all the $k$ smallest eigenvalues of the voltage-based Laplacian matrix are smaller than the $(k+1)$-th smallest eigenvalue of the ideal Laplacian matrix.} Eq.~\ref{Eq:Regular_Spectral_Clustering} represents the single-view spectral clustering cost minimization problem. In Eq.~\ref{Eq_Single_voltage_view}, superscript $(v)$ represents voltage.
\begin{align}
&\min_{{H^{(v)}} \in \mathbb{R}^{N \times k} } \Tr \left ( {H^{(v)}}^T \mathcal{L}^{(v)} {H^{(v)}} \right ) \label{Eq_Single_voltage_view} \\
& \text{subject to ${H^{(v)}}^T {H^{(v)}} = I$}. \nonumber
\end{align}
\textcolor{black}{The} assumption \textcolor{black}{may} only be partially relaxed} since location information does not give the complete information on meter-transformer \textcolor{black}{mapping}, and voltage information is also needed. Thus, geographical information complements the information from voltages to improve meter-transformer mapping. 


\textcolor{black}{The problem in Eq.~\ref{Eq_Single_voltage_view}} is minimized when the columns of $H^{(v)}$ are the eigenvectors of the Laplacian $\mathcal{L}^{(v)}$. In order to integrate multiple views, \textcolor{black}{the output of both views needs to be the same.} Therefore, we propose to devise an objective to minimize the disagreement between the two $H$'s. Let $H^{(l)}$ represent the GPS location information. Similar to finding the connectivity information in $V$ by constructing the similarity matrix $M$, the connectivity information in $H$ is obtained by constructing a similarity matrix  $K_{H}$. We choose a linear kernel $k(\mathbf{h}_i,\mathbf{h}_j) = \mathbf{h}_i^{T}\mathbf{h}_j,$  and so, $K_{H} = H^{(v)}{{H}^{(v)}}^{T}$ \textcolor{black}{is constructed}. The reason for choosing a linear kernel to measure the similarity of $H^{(\cdot)}$ is that the measure used in the similarity matrix for spectral clustering already  takes care of the non-linearities present in the data. 

Let us consider an ideal condition where voltage and location provide identical connectivity information. Therefore, $K_{{H}^{(v)}}$ and $K_{{H}^{(l)}}$
 would be similar, and the difference $K_{{H}^{(v)}} - K_{{H}^{(l)}}$ would be very small. \textcolor{black}{So}, to satisfy both conditions, we include  $K_{{H}^{(v)}} - K_{{H}^{(l)}}$ in the spectral clustering objective function with a suitable regularization coefficient. However, both  $K_{{H}^{(v)}}$ and $K_{{H}^{(l)}}$ \textcolor{black}{may} have different scales due to different data sources, and therefore, \textcolor{black}{each is normalized.} \vspace{-3mm}

\begin{align}
D\left ( {H}^{(v)},{H}^{(l)} \right ) = \left \| \frac{K_{{H}^{(v)}}}{\left \| K_{{H}^{(v)}} \right \|^2_{F}} - \frac{K_{{H}^{(l)}}}{\left \| K_{{H}^{(l)}} \right \|^2_{F}} \right \|^2_{F},
\label{Eq:coregularization}
\end{align}where $K_{{H}^{(v)}}$ is the similarity matrix for ${H}^{(v)},$ and $\left \| \cdot \right \|^2_{F}$ denotes the Frobenius norm of the matrix $H$, which is similar to the usual Euclidean norm by treating the matrix as a vector. The similarity matrices are normalized using their Frobenius norms to make them comparable across views. Frobenius norm has a useful property $||A||^2_F = \Tr(A^T A) = \Tr(AA^T).$
By substituting the last property and the distributive property of $\Tr$  in Eq.~\ref{Eq:coregularization} and ignoring the additive and multiplicative terms that depend on ${H}^{(v)}$ or ${H}^{(l)}$ individually, \textcolor{black}{we obtain} \vspace{-3mm}

\begin{align*}
D\left ( {H}^{(v)},{H}^{(l)} \right ) &= H -\Tr \left ( K_{{H}^{(v)}}{K_{{H}^{(v)}}}^T \right ). \\
&= -\Tr \left(  {H}^{(v)}{{H}^{(v)}}^{T} {H}^{(l)}{{H}^{(l)}}^{T}  \right).
\end{align*}
\textcolor{black}{The above disagreement needs to be reduced} between the clusterings of voltage view $v$ and location view $l$. Combining the objectives of individual views with the above term via a \textcolor{black}{multiplier} $\lambda$, we \textcolor{black}{obtain} a joint \textit{minimization} problem for the two graphs, 
where the hyperparameter $\lambda$ trades off the spectral clustering objectives and the spectral embedding disagreement term. The joint optimization problem given by the above equation can be solved using alternating minimization with respect to ${H}^{(v)}$ and ${H}^{(l)}$. For a given ${H}^{(l)}$, we \textcolor{black}{get} the following optimization problem in ${H}^{(v)}$
\begin{align}
&\min_{{H}^{(v)} \in \mathcal{R}^{n \times k}} \Tr \left \{ {{H}^{(v)}}^{T} \left( \mathcal{L}^{(v)} - \lambda {H}^{(l)}{{H}^{(l)}}^{T}  \right) {H}^{(v)} \right \} \label{Eq:combined_laplacian} \\
&\text{subject to   }{H}^{(v)}{{H}^{(v)}}^{T} = I. \nonumber
\end{align}

\textcolor{black}{By comparing Eq.~\ref{Eq_Single_voltage_view} with Eq.~\ref{Eq:combined_laplacian}, we observe that Eq.~\ref{Eq:combined_laplacian} is} the regular spectral clustering objective function on view $v$ with graph Laplacian $\mathcal{L}^{(v)} - \lambda {H}^{(l)}{{H}^{(l)}}^{T}$. \textcolor{black}{It} can be seen as a way of combining kernels or Laplacians. The difference from standard kernel combination (kernel addition, for example) is that the combination is updated at each step, as guided by the clustering algorithm \cite{kumar2011multiview}. Using such a framework, we effectively combine the voltage and the location information for \textcolor{black}{meter-transformer mapping} identification with obtaining a common solution. 
\textcolor{black}{Below, two metrics are shown \textcolor{black}{to compute} geographical distance, and the best metric is selected.} 
\vspace{-5mm}

\subsection{Metric Evaluation for Geographical Distance}
As discussed, location is important because electric lines underground and overhead follow streets. The purpose is to use the GPS coordinates (latitudes and longitudes) $\mathbf{l}^1,\cdots,\mathbf{l}^N \in \mathbb{R}^{2 \times 1}$ to obtain the distance between meters.
\subsubsection{Euclidean distance-based metric}
For two points \textcolor{black}{$i$ and $j$}, using the difference of the latitudes and the longitudes, one can estimate the angle \textcolor{black}{$\theta_{ij}$} between the points (subtended at the center of the Earth).
\begin{align}
    d_{ij} = R_E \times \theta_{ij} := R_E \times \sqrt{\left (l^i_1 - l^j_1 \right)^2 + \left (l^i_2 - l^j_2 \right)^2},
    \label{Eq:Euclidean_Angle}
\end{align}
where $d_{ij}$ is the distance between the $i$-th and the $j$-th points, $R_E$ is the radius of Earth, $l_1$ and $l_2$ are the latitude and longitude in radians, respectively.

\textcolor{black}{\textit{Drawback:} Although the approach is easy to compute, the distance is \textcolor{black}{inaccurate.} For example, the latitudes and longitudes have the same weightage, which is true near the equator. \textcolor{black}{However, near the poles, the same change in longitude has a much lesser effect than an equal latitude change.} Therefore, Eq.~\ref{Eq:Euclidean_Angle} is not valid for calculating the distance.}\vspace{-3mm}

\textcolor{black}{\subsubsection{Haversine Distance} To obtain the exact distance, consider two points with the same latitude $l^i_1 = l^j_1 = l_1$. In this case, the distance between the two points will be $d_ij = R_E \times \left ( l^i_2 - l^j_2 \right ) \times \cos{l_1}$. However, when the latitude and longitude both change, we obtain the distance using the Haversine formula. Let $A_1 := \sin^2 \left( \frac{ \Delta l_1 }{2} \right)$ and $A_2 := \sin^2 \left( \frac{ \Delta l_2 }{2} \right)$.
\begin{align*}
    d_{ij} = 2 \times R_E \times
    \arcsin \left ( \sqrt{ A_1 +\cos(l^i_1) \cos(l^j_1) A_2 } \right ) 
\end{align*}
} 
\textit{Remark. The spectral embedding-based meter-transformer (M.T.) mapping identification method is guaranteed to recover M.T. mapping under the assumption on the Laplacian matrix, i.e., \textcolor{black}{all the $k$ smallest eigenvalues of the voltage-based Laplacian matrix are smaller than the $(k+1)$-th smallest eigenvalue of the ideal Laplacian matrix. Moreover, we partially relax the assumption by using meter location information to aid meter voltage information via multi-view spectral clustering.}}



\vspace{-3mm}
\section{Numerical Results}

\begin{figure*}[ht]
\centering
\subfloat[Subfigure $1$ list of figures text][Result of the BIRCH clustering algorithm. \textcolor{black}{BIRCH algorithm identifies three clusters instead of two clusters as it is not suitable for voltage data since the radius hyperparameter needs to be hardcoded. 
}]{
\includegraphics[width=0.3\textwidth]{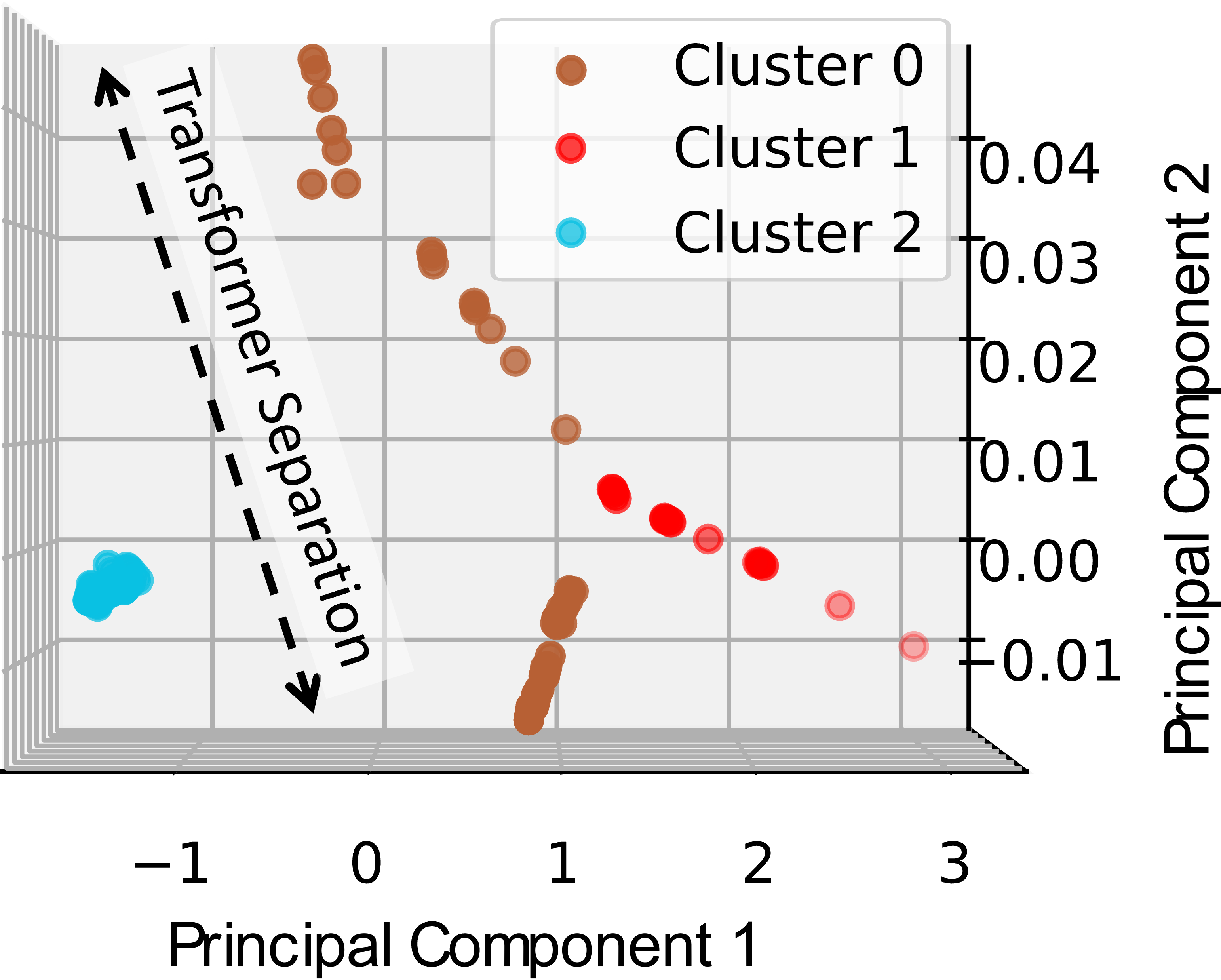}
\label{fig:Birch_Combined_IEEE123}}
\quad
\subfloat[Subfigure 1 list of figures text][Result of the $k-$means clustering algorithm. We specify the number of clusters in $k-$means, \textcolor{black}{but it identified incorrect clusters since the within-cluster optimization approach is unsuitable for voltage data.}]{
\includegraphics[width = 0.3\textwidth]{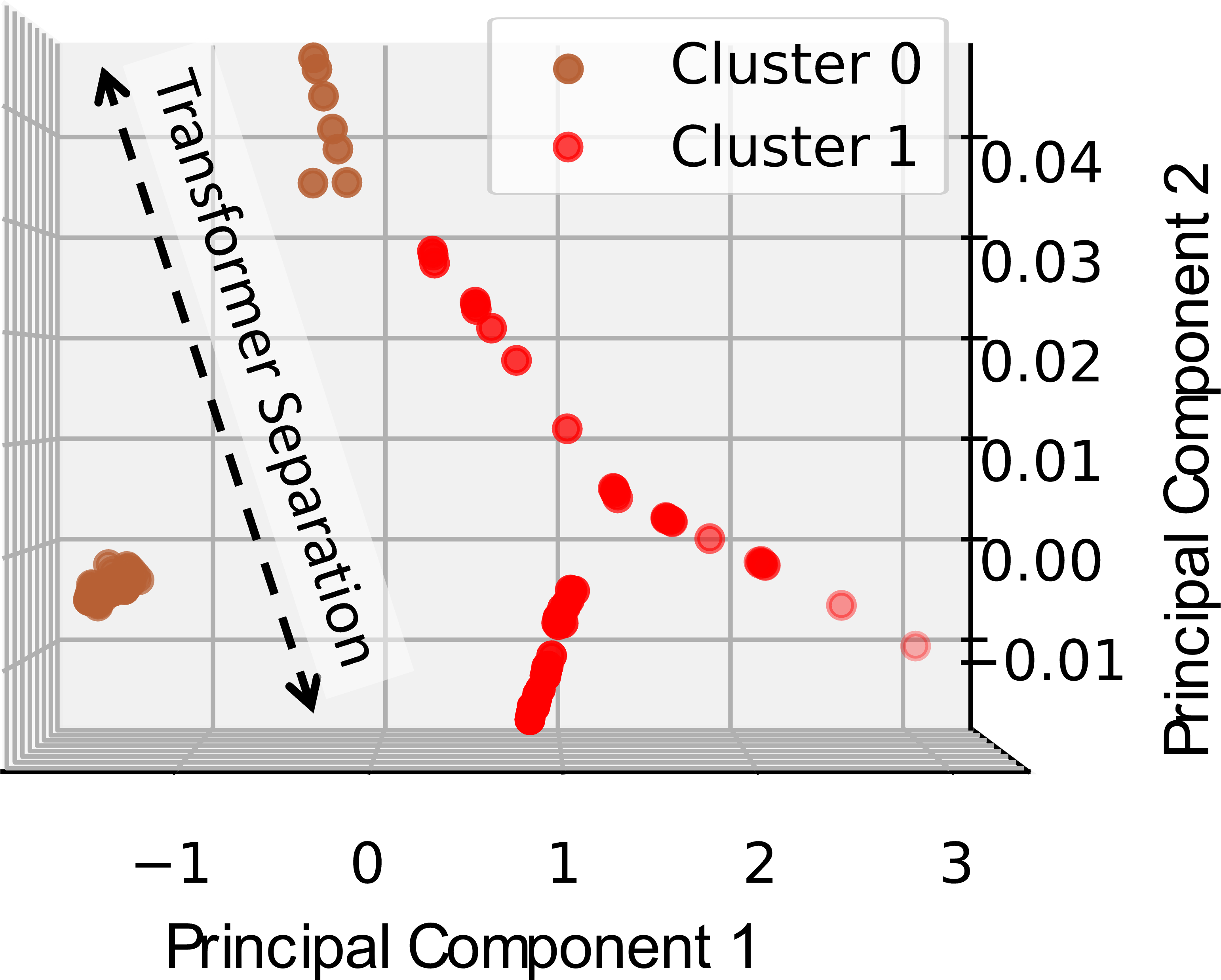}
\label{fig:Kmeans_Combined_IEEE123}}
\quad
\subfloat[Subfigure 2 list of figures text][Result of spectral clustering algorithm with voltage mutual information. Spectral clustering gave correct results for the transformer separation \textcolor{black}{using voltage data alone.}]{
\includegraphics[width=0.3\textwidth]{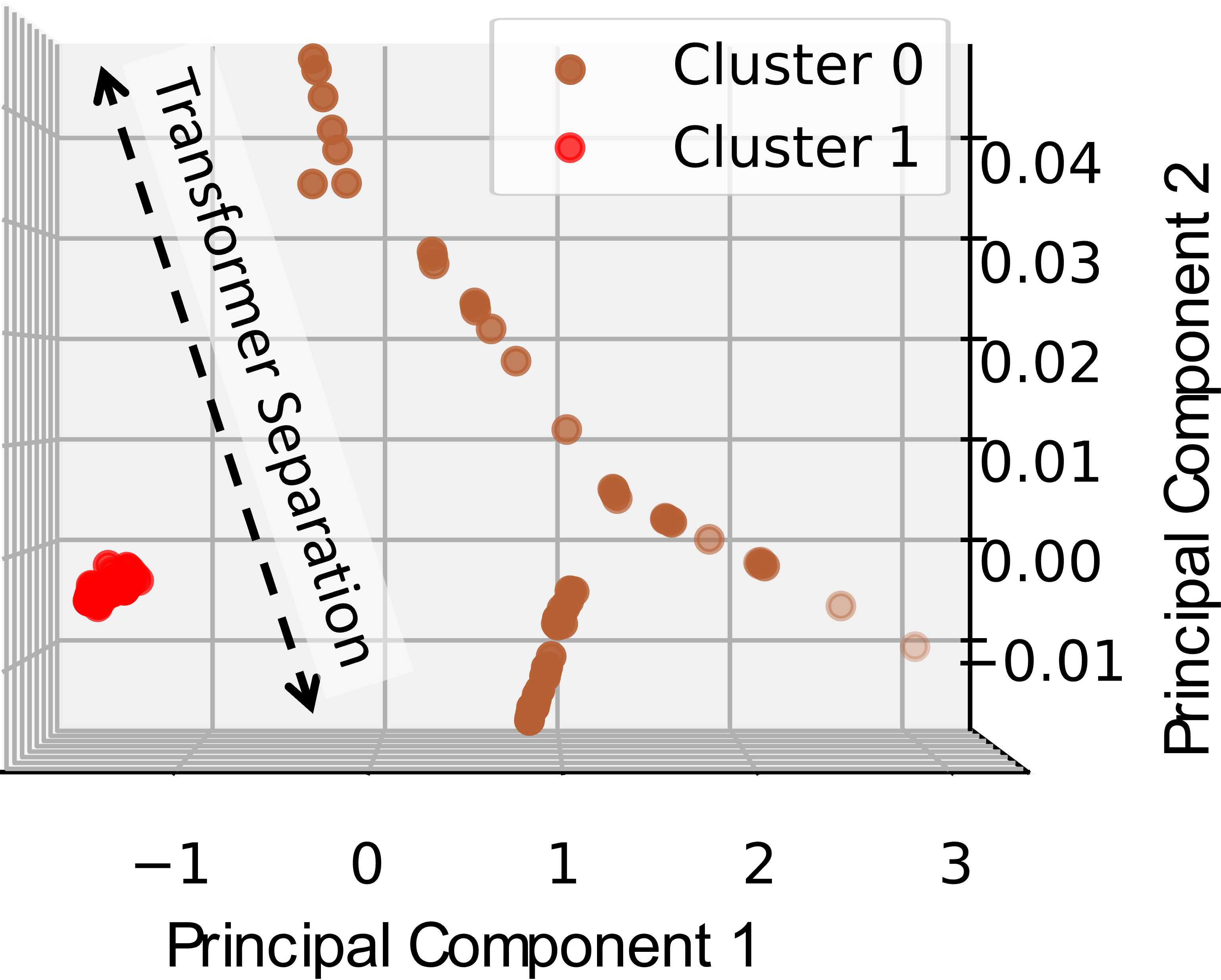}
\label{fig:SpectralClustering_Combined_IEEE123}}

\caption{Comparison of the three clustering algorithms using only voltage data on IEEE-$123$ bus test feeder \textcolor{black}{divided into two parts, as shown in Fig.~\ref{fig:IEEE-123_two_segments}. The dotted line for ``Transformer Separation'' shows the ground truth of the two clusters.}}
\label{fig:AlgoComparison.png}
\vspace{-4mm}
\end{figure*}

In this section, the proposed spectral clustering for meter-transformer mapping \textcolor{black}{is extensively validated using various test cases and real-world scenarios. Moreover, we also validate the assumption needed for spectral clustering on a real dataset.} To numerically validate the performance, we compare spectral clustering with representative clustering algorithms having the same setup as our proposed method: $k$-means (representing cost minimization methods) 
and BIRCH (representing \textcolor{black}{fixed-cluster size} methods), as mentioned in Section~\ref{sec:DBSCAN_Kmeans_BIRCH_SpecClust}. To compare performance, \textcolor{black}{the standard IEEE-$123$ bus system is modified} into various clusters by randomly dividing \textcolor{black}{so} that the structure remains preserved. Since \textcolor{black}{the goal is to have a method with a theoretical assurance, our guarantee is validated under the assumption mentioned in} Section~\ref{Sec:Spectral_Cluster_Guarantee}. \textcolor{black}{In addition,} \textcolor{black}{the} partial relaxation of the assumption involving transformer impedance and net consumptions, shown in Section~\ref{Sec:Co-regularized_multi-view_Spectral_Clustering}, \textcolor{black}{is also validated}. \textcolor{black}{We use the IEEE-123 bus system and real distribution systems from our partner utility for validation.}
\vspace{-2mm}
\subsection{Data Description and Algorithm Validation Strategy}
\label{Subsec_Data_Description_Algorithm_Validation_Strategy}

To make the simulation on test cases closer to reality, historical AMI consumption data from a local utility is used. Having the data and the customized test systems, we used power flow to generate voltage data from historical consumption data using OpenDSS for our algorithm validation. The algorithm is implemented using MATPOWER and Scikit-learn libraries.

In addition to the generated voltage data, \textcolor{black}{the proposed algorithm is validated} on historical voltage data from our partner utility. Set one is voltage data for Feeder A for the entire of $2018$.  Set two is voltage data of all other feeders for December $2018$. Time resolution is $15$ minutes, and the data availability is shown in  Table~\ref{tab:utility_data}. \textcolor{black}{Algorithmic meter-transformer mapping is validated using the ground truth of the same area obtained from the utility.}

\begin{table}[htbp]
\vspace{-4mm}
  \caption{Voltage dataset 1 and dataset 2 from our partner utility with a 15-minute resolution.}
  \resizebox{.47\textwidth}{!}{
  \centering
    \begin{tabular}{lcc}
    \toprule
    \textbf{Feeder A} & \textbf{$\#$ of meters} & \textbf{Total data number }  \\
    \midrule
    \textbf{January } & $153$ & $153 \times 2,976=455,328$\\
    \textbf{February } & $165$ & $165 \times 2,688=443,520$ \\
    \textbf{March } & $171$ & $171 \times 2,976=508,896$\\
    \textbf{April } & $173$ & $173 \times 2,880=498,240$ \\
    \textbf{May } & $195$ & $195 \times 2,976=580,320$ \\
    \textbf{June } & $199$ & $199 \times 2,880=573,120$ \\
    \textbf{July } & $207$ & $207 \times 2,976=616,032$ \\
    \textbf{August } & $209$ & $209 \times 2,976=621,984$ \\
    \textbf{September } & $210$ & $210 \times 2,880=604,800$ \\
    \textbf{October } & $212$ & $212 \times 2,976=630,912$ \\
    \textbf{November } & $213$ & $213 \times 2,880=613,440$ \\
    \textbf{December } & $213$ & $213 \times 2,976=633,888$ \\
    \midrule
    \textbf{All Other Feeders } & $81097$ & $81,097 \times 2,976=241,344,672$ \\
    \bottomrule
    \end{tabular}%
    }
  \label{tab:utility_data}%
   \vspace{-2mm}
\end{table}%
In addition, a third dataset containing five days of smart meter voltage data from our partner utility \textcolor{black}{is received. It contains data} from 1st September 2019 to 5th September 2019. A total of $5593600$ smart meter IDs are present in the data. 
About $1.73\%$ of the data has missing values. In addition to the smart meter data, the locations of $163518$ distribution transformers \textcolor{black}{are also received}.
\vspace{-2mm}

\subsection{\textcolor{black}{Validation of the Assumption}}
\label{Sec_Validation_of_the_Three_Assumptions}

\begin{figure}[ht]
    \centering
    \includegraphics[width = 0.45 \textwidth]{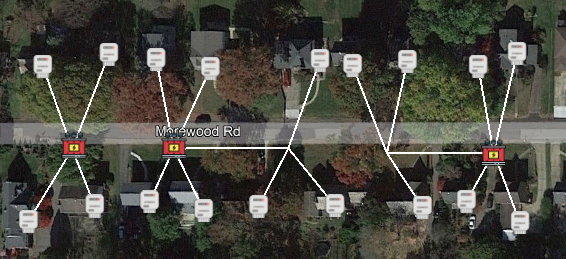}
    \vspace{-3mm}
    \caption{\textcolor{black}{shows a test case with nearby transformers for validating the assumption.}}
    \label{fig:assumption_validation_testcase}
\end{figure}

\textcolor{black}{Before we validate the algorithm, we also need to validate the assumption we made for the guarantee of the algorithm. The assumption is that the consumptions of inter-transformer meters are not the same, and the noise in voltage data is limited so that all the $k$ smallest eigenvalues of the voltage-based Laplacian matrix are smaller than the next smallest eigenvalue of the ideal Laplacian matrix. 
To validate the assumption, we consider real scenarios with nearby transformers. For example, Fig.~\ref{fig:assumption_validation_testcase} shows a test case with nearby transformers.}

\textcolor{black}{
The guarantee exists if the $k$ smallest eigenvalues of the voltage-based Laplacian matrix are smaller than the next smallest eigenvalue of the ideal Laplacian matrix. Validating the assumption needs a real Laplacian matrix and an ideal Laplacian matrix. The real Laplacian matrix \textcolor{black}{is} constructed from the voltage data via the similarity matrix. The ideal Laplacian matrix \textcolor{black}{is} formed from the unweighted graph adjacency matrix using ground truth information. For example, an unweighted graph adjacency matrix is a binary square matrix, with $1$ for meters supplied by the same transformer and $0$ otherwise.}

\begin{figure}[ht]
    \centering
    \includegraphics[width = 0.45 \textwidth]{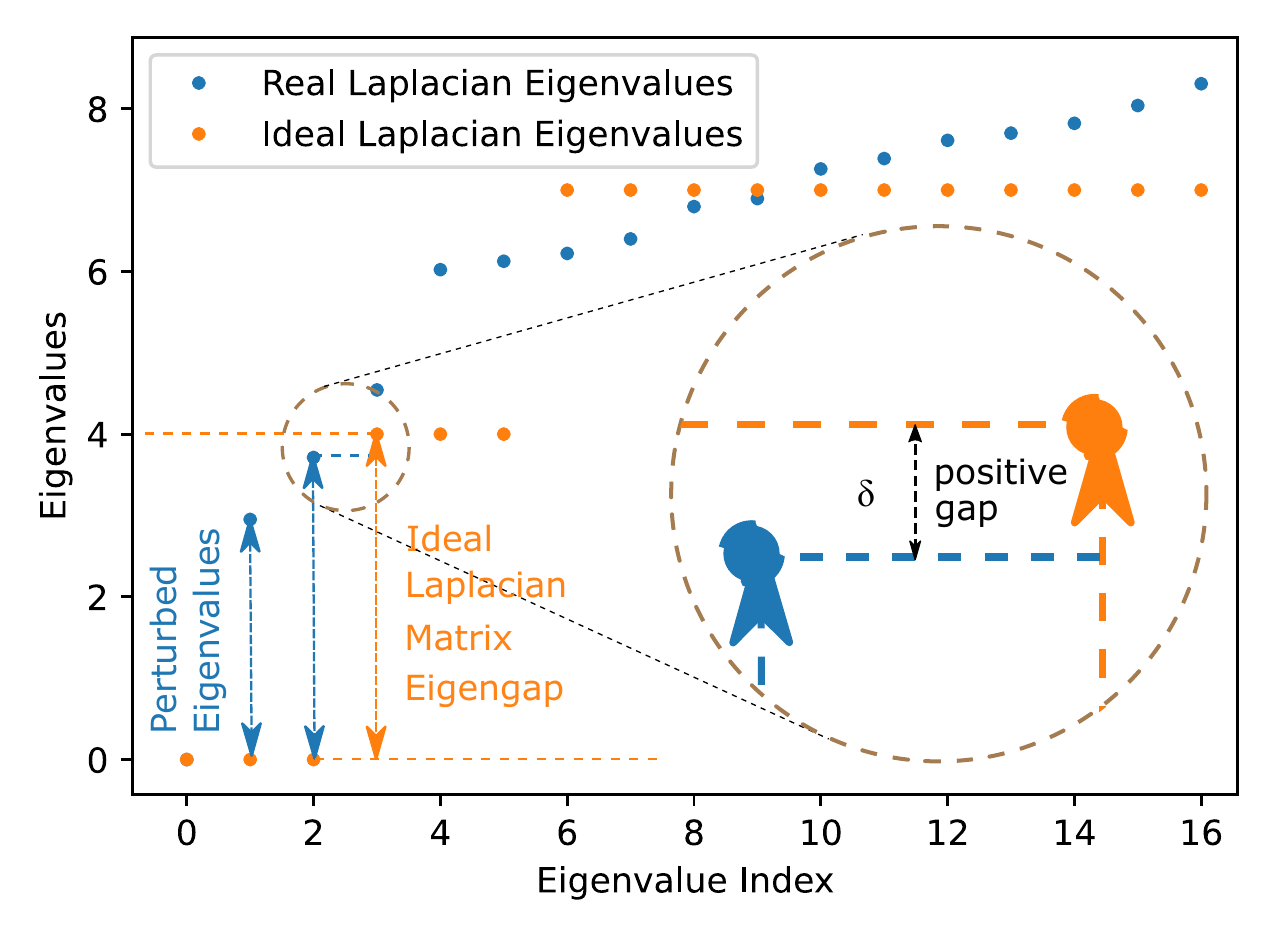}
    \vspace{-3mm}
    \caption{\textcolor{black}{shows a scatter plot of eigenvalues of ideal and real Laplacian matrices. We observe positive $\delta$ and validate the assumption for the test case in Fig.~\ref{fig:assumption_validation_testcase}.}}
    \label{fig:validation_assumption2.png}
\end{figure}

\textcolor{black}{Fig.~\ref{fig:validation_assumption2.png} demonstrates the validation of the assumption for the test case of Fig.~\ref{fig:assumption_validation_testcase}. For the ideal Laplacian matrix, we observe the first $3$ eigenvalues equalling zero, which represents $3$ clusters. Moreover, we observe the eigengap equal to $4$. For the real Laplacian matrix, we observe that the second and third eigenvalues are non-zero. Moreover, we observe that the gap $\delta$ between the third eigenvalue of the ideal Laplacian matrix and the fourth eigenvalue of the real Laplacian matrix is positive. Therefore, the assumption is valid.}

\begin{figure}[ht]
    \centering
    \includegraphics[width = 0.45 \textwidth]{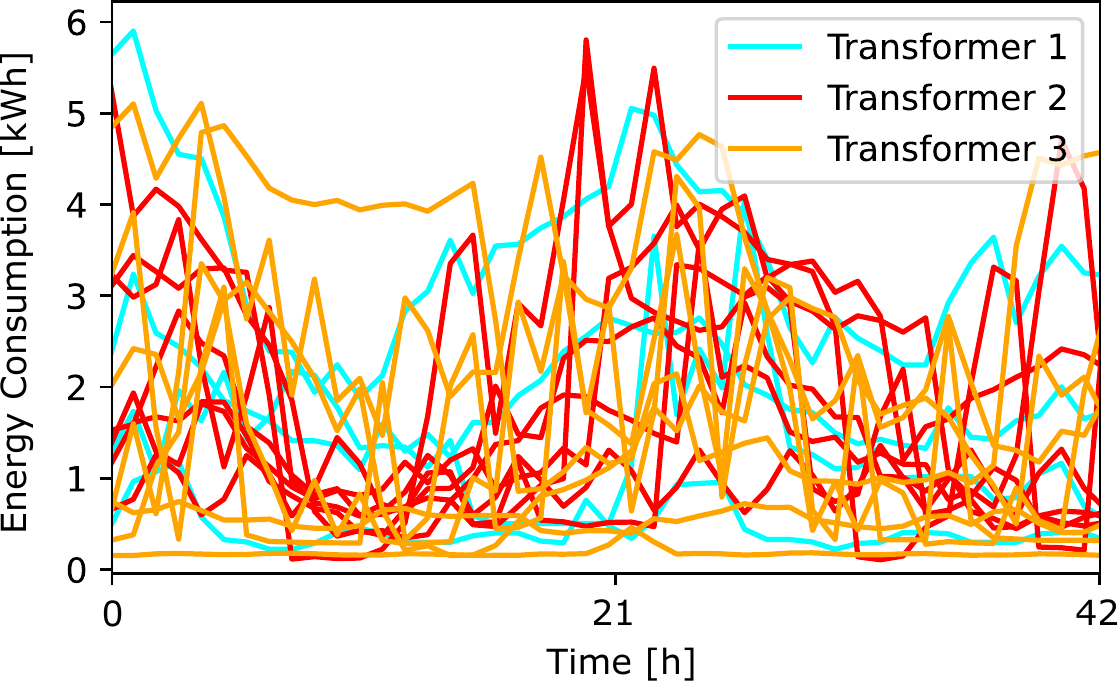}
    \vspace{-3mm}
    \caption{\textcolor{black}{shows a line plot of meter energy consumptions for the test case of Fig.~\ref{fig:assumption_validation_testcase}.}}
    \label{fig:validation_assumption3}
\end{figure}

\begin{figure}[ht]
    \centering
    \includegraphics[width = 0.45 \textwidth]{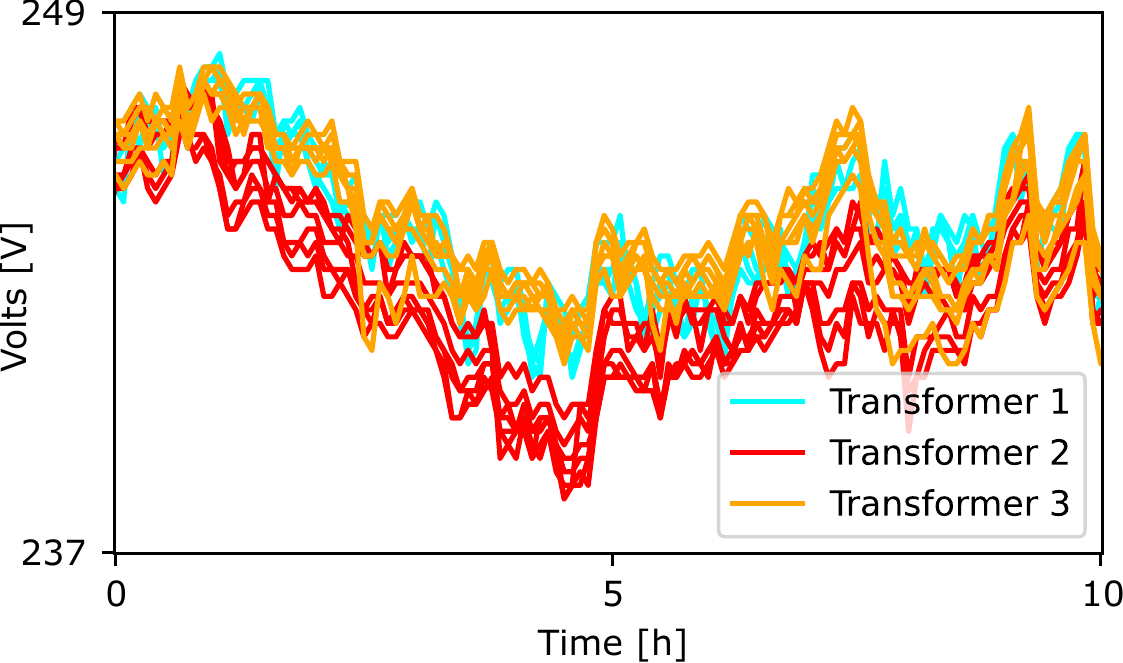}
    \vspace{-3mm}
    \caption{\textcolor{black}{shows a line plot of meter voltages for the test case of Fig.~\ref{fig:assumption_validation_testcase}.}}
    \label{fig:validation_assumption1.png}
\end{figure}

\textcolor{black}{The valid assumption implies that the consumptions of inter-transformer meters are not the same, as observed in Fig.~\ref{fig:validation_assumption3}.}
\textcolor{black}{Moreover, the meter voltages are similar for meters supplied by the same transformer and otherwise different, as seen in Fig.~\ref{fig:validation_assumption1.png}. The assumption for the guarantee of M.T. identification using spectral clustering is validated. Next, we show validation of the proposed method.}

\subsection{Validation of Robustness for Spectral Clustering}
\label{Subsec_Validation_of_Robustness_for_Spectral_Clustering}
\textcolor{black}{In the following, we validate spectral clustering using IEEE test systems and real systems from our partner utility.}

\subsubsection{Validation \textcolor{black}{using} the IEEE-123 Test Case}
IEEE-123 bus system operates at $4.16kV$. It contains overhead and underground lines, switches, shunt capacitor banks, unbalanced constant current, constant power, and impedance loading~\cite{Schneider2018}.

\begin{figure}[htbp]
    \centering
    \includegraphics[width = 0.4 \textwidth]{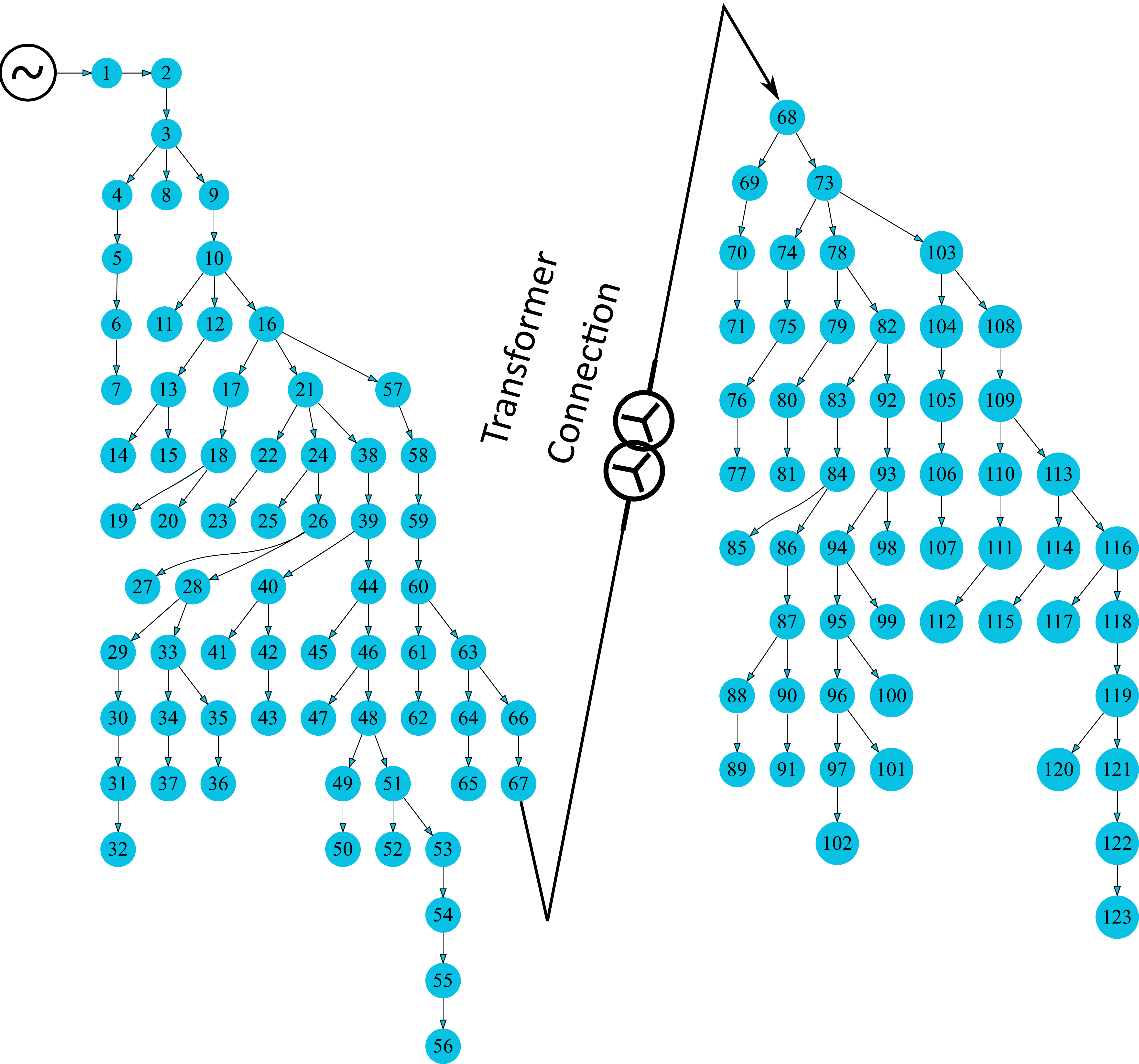}
    \vspace{-2mm}    
    \caption{shows the IEEE-123 bus test feeder divided into two parts by a transformer between them.}
    \label{fig:IEEE-123_two_segments}
    \vspace{-2mm}
\end{figure}

For numerical validation of robustness, the proposed method \textcolor{black}{is validated} on the IEEE-$123$ bus case.
In particular, Fig.~\ref{fig:IEEE-123_two_segments} shows a transformer connected between buses $67$ and $68$. A voltage dataset is generated using load-flow analysis utilizing a real consumption dataset from our partner utility using the MATPOWER library. Fig.~\ref{fig:AlgoComparison.png} \textcolor{black}{compares} spectral clustering with two famous clustering approaches: $k-$means and BIRCH. We observe that only our proposed method, shown in Fig.~\ref{fig:SpectralClustering_Combined_IEEE123}, consistently outperforms both famous methods for any placement of the transformer and various voltage datasets using randomly selected datasets of the load profiles from our partner utility. \textcolor{black}{Such a superiority is due to the proposed method using the eigenvectors of the graph Laplacian} \textcolor{black}{to cluster high-dimensional smart meter datasets effectively.}

\subsubsection{Validation \textcolor{black}{using} Utility Systems}
\begin{figure}[htbp]
    \centering
    \includegraphics[width = 0.2 \textwidth]{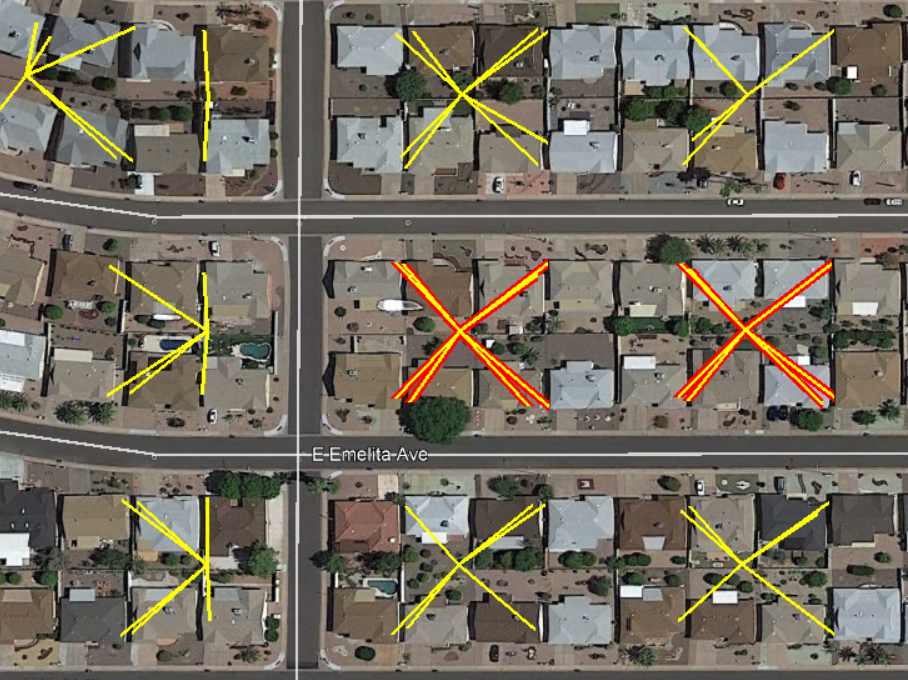}
    \vspace{-2mm}
    \caption{The obtained \textcolor{black}{meter-transformer mapping} (red color) coincides with the ground truth topology (yellow) for a utility test case.}
    \label{fig:Good_Example_1}
    \vspace{-4mm}
\end{figure}

In addition to the IEEE-123 bus test system, the proposed method \textcolor{black}{is also validated} on real utility \textcolor{black}{systems.} For example, Fig.~\ref{fig:Good_Example_1} shows the real utility system. \textcolor{black}{As} one can see that the results match the ground truth. \textcolor{black}{Moreover, the proposed method successfully recovers the ground truth for the test case of Fig.~\ref{fig_meter_near_wrong_cluster}. Therefore, spectral clustering is more robust than other methods. Next we empirically validate the robustness guarantee.}

\vspace{-3mm}
\subsection{Validation of Spectral Clustering Robustness Guarantee}
\label{Subsec_Validation_Spectral_Clustering_Robustness_Guarantee}

\begin{figure}[htbp]
    \centering
    \vspace{-3mm}
    \includegraphics[width = 0.3
    \textwidth]{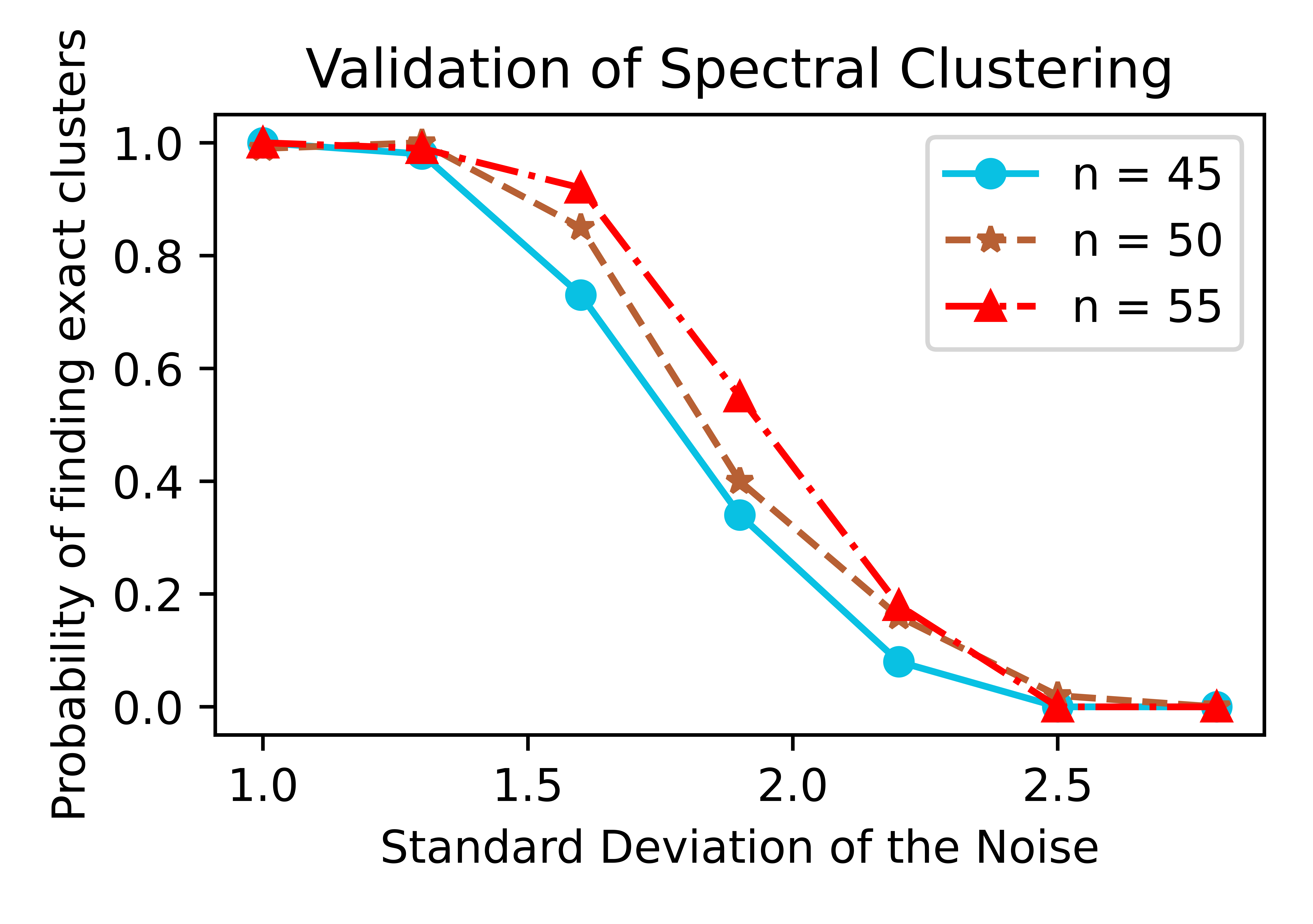}
    \vspace{-3mm}
    \caption{The \textcolor{black}{empirical} probability of success in obtaining the correct clusters reduces with the degree of noise penetration. 
    }
    \label{fig:Spectral_Clustering_Guarantee}
    \vspace{-2mm}
\end{figure}

Section~\ref{Sec:Spectral_Cluster_Guarantee} guarantees that the spectral clustering algorithm is robust 
under \textcolor{black}{the assumption mentioned.}
\textcolor{black}{The guarantee is numerically validated on the IEEE-$123$ bus system} shown in Fig.~\ref{fig:IEEE-123_two_segments}. After generating voltages using the system, $n$ \textcolor{black}{neighboring} bus voltages \textcolor{black}{are randomly selected} from the two parts of the IEEE-123 bus feeder, where $n$ is shown in Fig.~\ref{fig:Spectral_Clustering_Guarantee}. \textcolor{black}{Further, noise is added to \textcolor{black}{observe} the robustness of the algorithm.} \textcolor{black}{Finally,} the figure presents the \textcolor{black}{empirical} probability of success in finding the correct clusters. For example, the probability of finding the correct clusters is a monotonically decreasing function \textcolor{black}{of the added noise}. 
However, as one can see, below a certain noise level, \textcolor{black}{the guarantee of finding the \textcolor{black}{ground truth} clusters exists.}

\vspace{-3mm}







\subsection{Validation of Improvement by Co-regularized Multi-view Spectral Clustering with Distance}
\label{Subsec_Validation_Improvement_Coregularized_Multiview_Spectral_Clustering_with_Distance}
\textcolor{black}{To generalize the applicability of our algorithm, we partially relax the assumption, 
as mentioned in Section~\ref{Sec:Co-regularized_multi-view_Spectral_Clustering}. Specifically, we use meter location information to partially relax the condition\textcolor{black}{, i.e., all the $k$ smallest eigenvalues of the voltage-based Laplacian matrix are smaller than the $(k+1)$-th smallest eigenvalue of the ideal Laplacian matrix.}}

\subsubsection{Validation \textcolor{black}{using} the IEEE-123 Test Case}
\begin{figure}[htbp]
    \centering
    \includegraphics[width = 0.3 \textwidth]{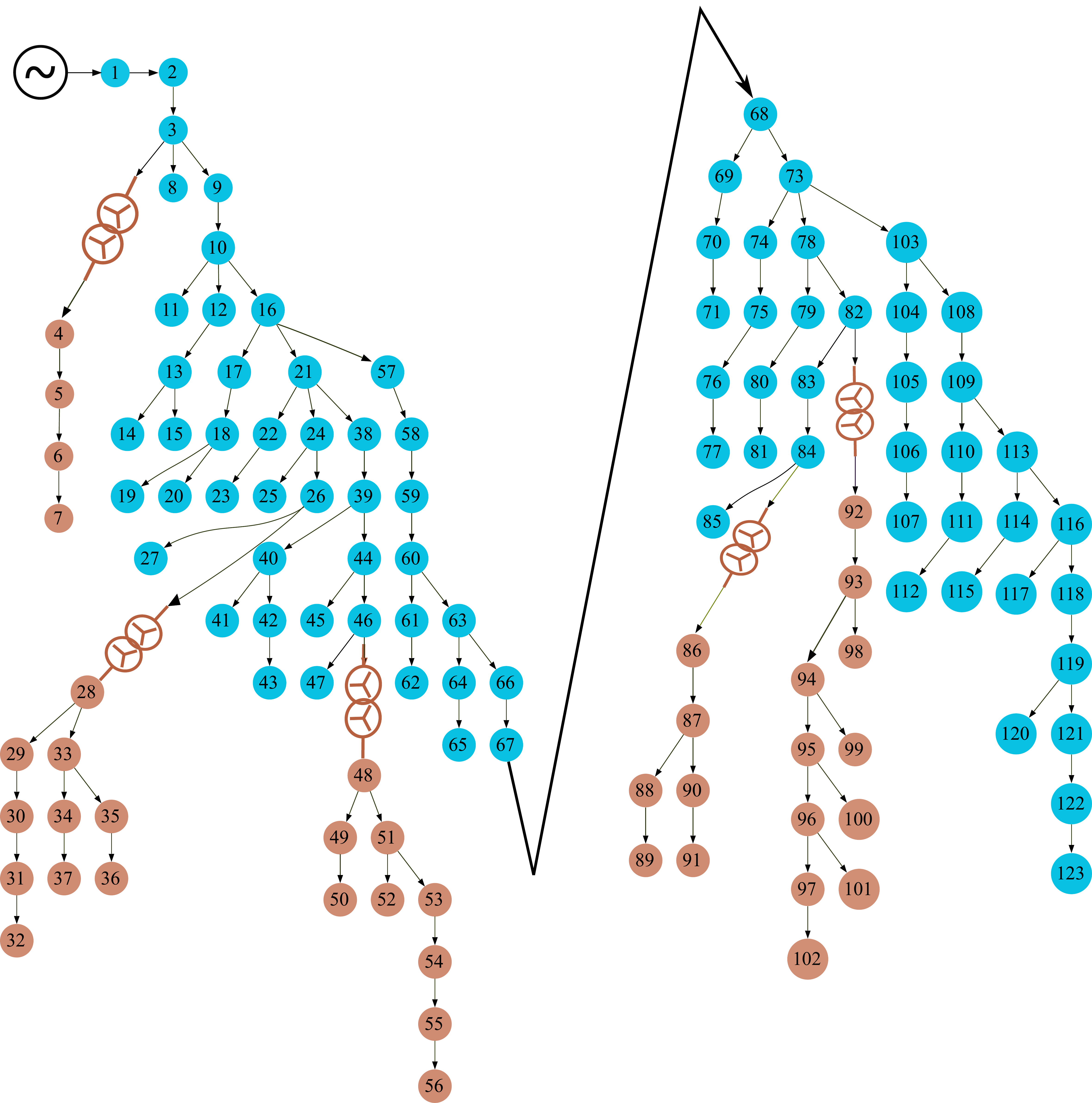}
    \caption{shows the same IEEE-123 bus test feeder with five transformers. The secondaries are shown in brown color.}
    \label{fig:IEEE-123_6_segments}
\end{figure}

\begin{figure}[htbp]
\vspace{-5mm}
\centering
\subfloat[Subfigure $1$ list of figures text][Regular spectral clustering \textcolor{black}{incorrectly clusters the highlighted data point.}]{
\includegraphics[width=0.2\textwidth]{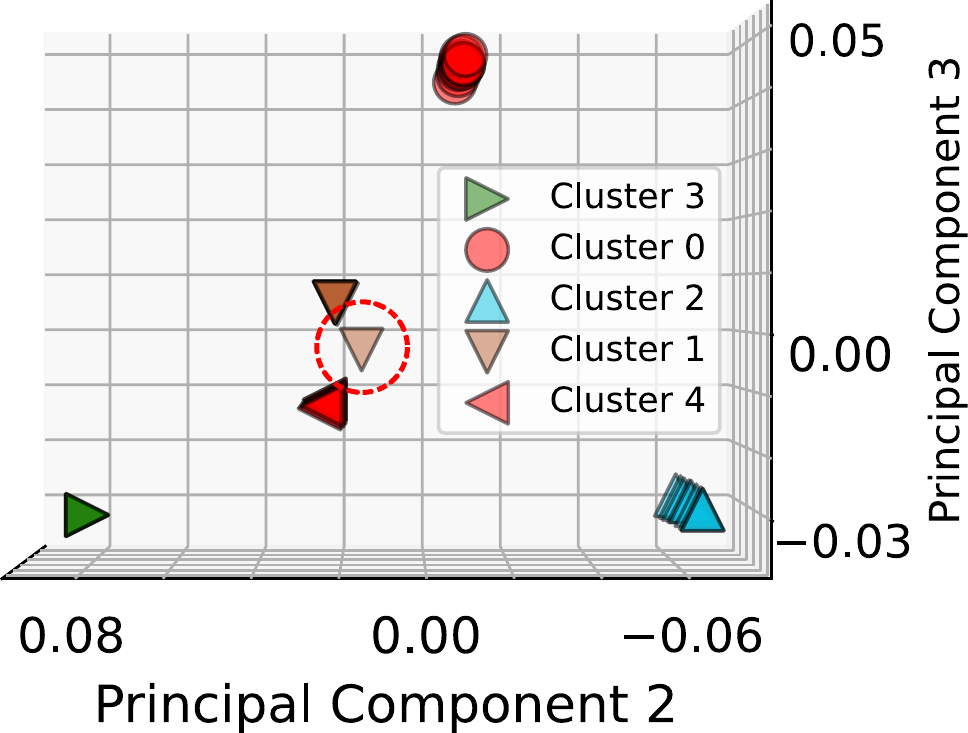}
\label{fig:Regular_Spectral_Clustering_IEEE-123_bus_VS_multiview}}
\quad
\subfloat[Subfigure 1 list of figures text][Multi-view spectral clustering \textcolor{black}{correctly clusters the highlighted data point.}]{
\includegraphics[width = 0.2\textwidth]{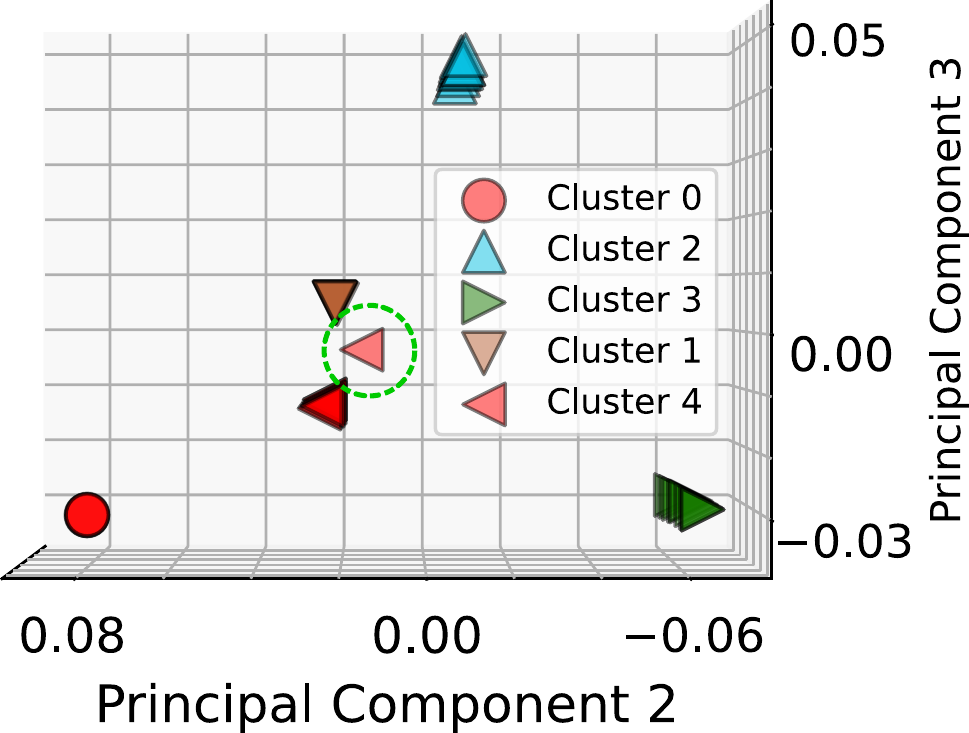}
\label{fig:Spectral_multiview_6_clusters}}
\caption{shows the comparison of the two algorithms on the IEEE-123 bus test feeder with five transformers, as shown in Fig.~\ref{fig:IEEE-123_6_segments}. Clusters are represented by colors and shapes.}
\label{fig:Multiview_Comparison_IEEE123}
\vspace{-2mm}
\end{figure}

To simulate such a scenario on the IEEE-123 test system, \textcolor{black}{we need a geographical view to improve voltage-based clustering. So, we assign} geographical coordinates (latitudes and longitudes) to each bus in the IEEE-123 bus test case. 
\textcolor{black}{Then,} \textcolor{black}{we insert transformers in the IEEE-123 bus case in such a way so that the resulting smart meter clusters are far from each other,} as shown in Fig.~\ref{fig:IEEE-123_6_segments}. The regular spectral clustering algorithm incorrectly resolves the clusters, as shown in Fig.~\ref{fig:Regular_Spectral_Clustering_IEEE-123_bus_VS_multiview}. \textcolor{black}{However, we obtain the correct results using multi-view spectral clustering,} as shown in Fig.~\ref{fig:Spectral_multiview_6_clusters}.

\textcolor{black}{Fig.~\ref{fig:IEEE-123_6_segments} shows a complex meter-transformer mapping scenario where the assumption that all the $k$ smallest eigenvalues of the voltage-based Laplacian matrix are smaller than the $(k+1)$-th smallest eigenvalue of the ideal Laplacian matrix 
does not satisfy. Therefore, regular spectral clustering is not guaranteed to work, and it does not work, as shown in Fig.~\ref{fig:Regular_Spectral_Clustering_IEEE-123_bus_VS_multiview}. Moreover, incorporating the location data via multi-view spectral clustering identifies the meter-transformer mapping \textcolor{black}{accurately,} as shown in Fig.~\ref{fig:Spectral_multiview_6_clusters}. So far, we noted an improvement in spectral clustering using IEEE-123 system. Next, we show the improvement on real system.}

\subsubsection{Validation \textcolor{black}{using} Utility System} 
\begin{figure}[htbp]
\vspace{-4mm}
\centering
\subfloat[Subfigure $1$ list of figures text][Regular spectral clustering]{
\includegraphics[width=0.2\textwidth]{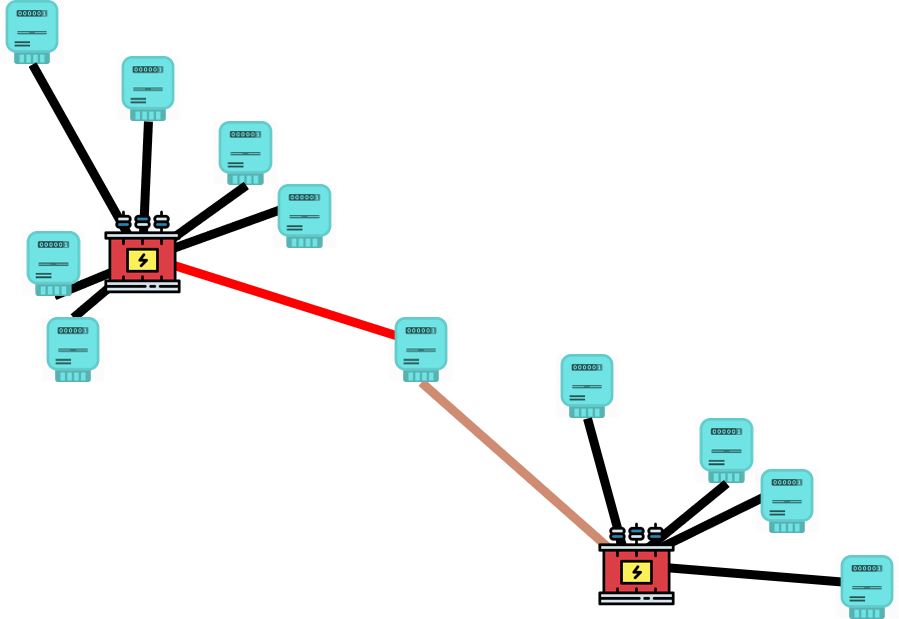}
\label{fig:APS_spectral_worse2}}
\quad
\subfloat[Subfigure 1 list of figures text][Multi-view spectral clustering]{
\includegraphics[width = 0.2\textwidth]{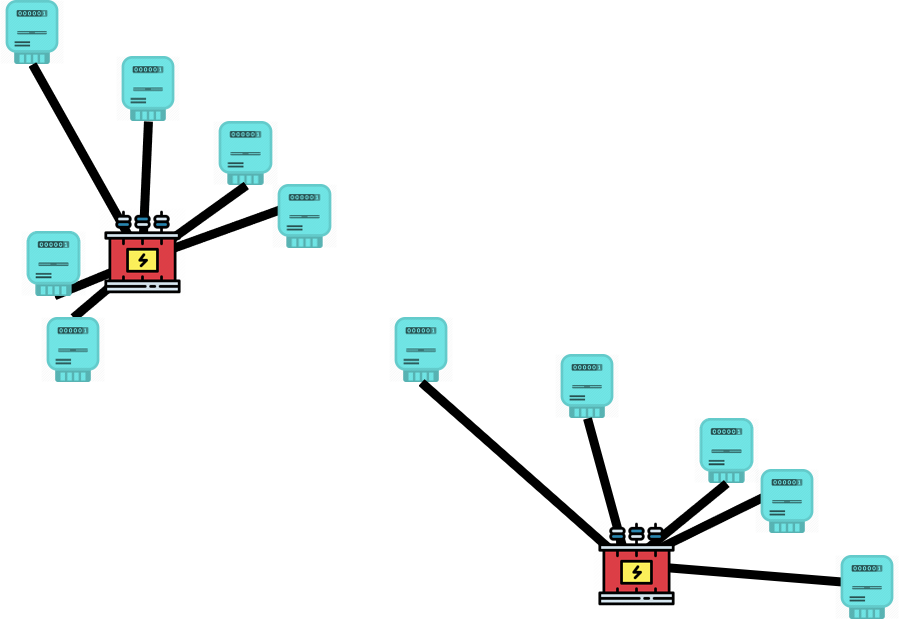}
\label{fig:APS_multiview_better}}
\caption{Comparison of regular spectral clustering versus co-regularized multi-view spectral clustering on a real utility system where the two transformers are far apart. Color code: red is incorrect prediction, yellow is missed prediction, and black is correct prediction.}
\label{fig:Multiview_Comparison_utility}
\end{figure}

Besides the IEEE test feeder, the co-regularized multi-view spectral clustering also improves the performance on real utility test cases, especially when the transformers are spaced apart. For example, Fig.~\ref{fig:APS_spectral_worse2} shows an example where the regular spectral clustering makes an incorrect connection. However, the correct results are obtained by incorporating the GIS information using co-regularized multi-view spectral clustering, as shown in Fig.~\ref{fig:APS_multiview_better}. \textcolor{black}{Therefore, the improvement is observed.} 


\textcolor{black}{
\textcolor{black}{To summarize, the validations in Sections  \ref{Sec_Validation_of_the_Three_Assumptions}, \ref{Subsec_Validation_of_Robustness_for_Spectral_Clustering}, and \ref{Subsec_Validation_Spectral_Clustering_Robustness_Guarantee}} \textcolor{black}{show} that this paper achieves a robust meter-transformer mapping with a robustness guarantee under the assumption mentioned. Moreover, the validations in Section~\ref{Subsec_Validation_Improvement_Coregularized_Multiview_Spectral_Clustering_with_Distance} demonstrates that the assumption on the Laplacian matrix spectrum is partially relaxed.}

\vspace{-2mm}
\section{Conclusion}




Distributed energy resources (DERs) have benefits, but they introduce challenges for the electric utilities, where real-time meter-transformer mapping can resolve them. \textcolor{black}{Past methods either try to identify topology with unreasonable assumptions or ignore common challenging scenarios, where a meter’s voltage may be similar to the non-parent transformer’s cluster. However, progressing from meter to meter can lead to the parent transformer’s cluster due to transformer impedances. Our proposed method utilizes such information via spectral embedding of voltage information. We also provide proof of a guarantee under a reasonable assumption on the voltage Laplacian matrix. Moreover, we employ meter location information to partially relax the assumption on the Laplacian matrix. Such a method is robust and easy to implement. Also, it does not assume a specific shape of transformer secondary circuits, a variety of data sensors, specific data distribution, or grid probing equipment.} The proposed method accurately identifies meter-transformer mapping on the IEEE test systems and real feeders from our partner utility.


\ifCLASSOPTIONcaptionsoff
  \newpage
\fi

\vspace{-3mm}
\bibliographystyle{IEEEtran}
\bibliography{Main.bib} 


\appendices
\section{Proof of Theorem 1}
\label{Appendix_Proof_of_Theorem_1}
The solution to the constraint minimization problem shown below \textcolor{black}{is the $H$ comprising $k$ eigenvectors corresponding to the $k$ smallest eigenvalues of $\mathcal{L}$.}
\begin{align*}
&\min_{H \in \mathbb{R}^{N \times k} } \Tr \left ( H^T \mathcal{L} H \right ) \\
& \text{subject to $H^T H = I$}. \nonumber
\end{align*}
\begin{IEEEproof}
\textcolor{black}{The Laplacian matrix $\mathcal{L}$ is a singular matrix with zero as an eigenvalue. Therefore,} the theorem is obvious if we consider all $\mathbf{h_i}$ as the eigenvectors of $\mathcal{L}$. \textcolor{black}{Generally, one can} resolve $\mathbf{h_i}$ into \textcolor{black}{components along the} eigenvectors of $\mathcal{L}$ \textcolor{black}{to prove the theorem.}
\end{IEEEproof}

\section{Proof of Theorem 4}
\label{Appendix_Proof_of_Theorem_4}

Let $X := \{\mathbf{x}_1, \mathbf{x}_2, \cdots , \mathbf{x}_k\} \in \mathbb{R}^{n \times k}$ be the matrix containing eigenvectors of $\mathcal{L}$ corresponding to eigenvalue $0$. 
\begin{enumerate}
\item The row vectors  of $X$ corresponding to the data points of the same cluster are equal. 
\item The row vectors of $X$ corresponding to different clusters are orthogonal.
\end{enumerate}


\begin{IEEEproof}
\textcolor{black}{The proof for 1) and 2) is as follows. The eigenvectors of the block-diagonal Laplacian $\mathcal{L}$ are a union of the appropriately zero-padded eigenvectors of individual diagonal blocks. Therefore, the row-vectors of $X$ corresponding to the data points of the same cluster are equal. Also, the row vectors of $X$ corresponding to different clusters are orthogonal ($\mathbf{x}_i^T \mathbf{x}_j = 0$). Moreover, due to the repeated eigenvalue $0$, $X$ can be replaced by $XQ$, where $Q$ is any $k \times k$ orthogonal matrix. Therefore, the row vectors of $X$ will be replaced by the row vectors of $Q$. Hence, 1) and 2) remain valid.}
\end{IEEEproof}

\section{Proof of Theorem 8}
\label{Appendix_Proof_of_Theorem_8}

Let $\mathcal{L}$ be symmetric. Also, let the columns of $X_1$ form an orthogonal basis for the simple invariant subspace of $\mathcal{L}$ for the eigenvalue $\lambda$. Moreover, let $\tilde{X}_1$ be the approximation of $X_1$, and $\tilde{P}_1 = \tilde{X}_1^T \mathcal{L} \tilde{X}_1$ become the approximation of $\lambda I$. \textcolor{black}{Moreover, let $\eig(\tilde{P}_1) \subseteq [a,b]$.} Let  $R = \mathcal{L}\tilde{X}_1 - \tilde{X}_1 \tilde{P}_1$ be the residual of the approximation $\tilde{X}_1$. Moreover, let 
\textcolor{black}{$sep \left([a,b],\eig(\Lambda_2)\right)>\delta$, where $sep(\cdot,\cdot)$ is the minimum distance over all elements of the two sets.} Then,
$$||\tan {\Theta} [\mathcal{R}(\tilde{X}_1),\mathcal{R}(X_1)] || \leq \frac{||R||}{\delta}.$$

\begin{IEEEproof}
Let us consider without loss of generality that $\Lambda$ and $ X_1^T \mathcal{L} X_1 = \lambda I $ have the same dimensions~\cite{stewart1990matrix}. 
Moreover, let us consider our frame of reference aligned with $\mathcal{R}(\tilde{X}_1)$. Therefore, we write
$\tilde{X}_1 = \binom{I}{0}.$

It is because $\tilde{X}_1$ should span the simple invariant subspace $\mathcal{R}(\tilde{X}_1)$ and not the orthogonal space of $\mathcal{R}(\tilde{X}_1)$. Let $\theta_1 \leq \cdots \leq \theta_n$ \textcolor{black}{be} the canonical angles between the original $\mathcal{R}(X_1)$ and the perturbed invariant subspace $\mathcal{R}(\tilde{X}_1)$. Therefore, in such a frame of reference, the column vectors of $X_1$ are at angles w.r.t. $\tilde{X}_1$. Moreover, the column vectors of $X_2$ are orthogonal to the respective column vectors of $X_1$. So, 
\begin{align*}
X =& 
\begin{pmatrix}
X_1 & X_2
\end{pmatrix} = 
\begin{pmatrix}
\Gamma  & -\Sigma \\ 
\Sigma  & \Gamma 
\end{pmatrix},
\end{align*}
where $\Gamma = \diag(\cos{\theta_i})$ consist of cosines in ascending order, $\Sigma = \diag(\sin{\theta_i})$ consists of sines in descending order. In the coordinate system of $\mathcal{R}(\tilde{X}_1)$, the Laplacian matrix can be generally partitioned as
\begin{align*}
\mathcal{L} = \begin{pmatrix}
\mathcal{L}_{11} & \mathcal{L}_{12}\\ 
\mathcal{L}_{21} & \mathcal{L}_{22}
\end{pmatrix}
\end{align*}
Moreover, in this coordinate system, the approximation of $\lambda I$ is given as $\tilde{P}_1 = \tilde{X}_1^T \mathcal{L} \tilde{X}_1 = \mathcal{L}_{11}$, and the resulting residual $R$ are given as
\begin{align}
R = & 
\begin{pmatrix}
\mathcal{L}_{11} & \mathcal{L}_{12}\\ 
\mathcal{L}_{21} & \mathcal{L}_{22}
\end{pmatrix}
\binom{I}{0} -  \binom{I}{0} \mathcal{L}_{11} = \binom{0}{\mathcal{L}_{21}}
\label{Eq_Equation_Residual}
\end{align}

Since $X_2^T \mathcal{L} = \Lambda_2 X_2^T$, it implies
$\begin{pmatrix}
-\Sigma & \Gamma
\end{pmatrix} \mathcal{L} = \Lambda_2 \begin{pmatrix}
-\Sigma & \Gamma
\end{pmatrix}.$
Therefore, left multiplying Eq.~\ref{Eq_Equation_Residual} by $X^T_2 = \begin{pmatrix}
-\Sigma & \Gamma
\end{pmatrix}$ and $X^T_1$, we respectively obtain
\begin{align}
\Gamma \mathcal{L}_{21} = \Sigma \mathcal{L}_{11} - \Lambda_2 \Sigma .
\label{Eq_Residual_Z_X2}\\
\Sigma \mathcal{L}_{21} = \lambda \Gamma - \Gamma \mathcal{L}_{11}.
\label{Eq_Residual_Z_X1}
\end{align}
Solving Eq.~\ref{Eq_Residual_Z_X2} and Eq.~\ref{Eq_Residual_Z_X1} simultaneously for $\mathcal{L}_{11}$ and using $\Gamma^2 + \Sigma^2 = I$, we get
\begin{align}
\mathcal{L}_{11} = \lambda \Gamma^2 + \Sigma \Lambda_2 \Sigma
\label{Eq_L_11}
\end{align}
From Eq.~\ref{Eq_L_11}, we observe that $\mathcal{L}_{11}$ is a diagonal submatrix since $\Gamma, \Sigma, \Lambda_2$ are diagonal submatrices. Therefore, the diagonal values of $\diag(\mathcal{L}_{11})$ are the eigenvalues $\eig(\mathcal{L}_{11})$. Consequently, $\diag(\mathcal{L}_{11}) = \eig(\tilde{P}_1) \subseteq [a,b]$.

The i-th diagonal element of Eq.~\ref{Eq_Residual_Z_X2} is given as
$\cos{\theta_i} l_{ii}^{(21)} = \sin{\theta_i} ( l_{ii}^{(11)} - \lambda_{ii}^{(22)} ) \geq \delta \sin{\theta_i}$. The matrix $\Lambda_2$ contains all eigenvalues of $\mathcal{L}$, except $\lambda$. Specifically, $\lambda_{ii}$ is an i-th diagonal value of $\Lambda_2$. Moreover, we have $( l_{ii}^{(11)} - \lambda_{ii}^{(22)} ) \geq \delta$ since $l_{ii}^{(11)} \in [a,b]$, and $sep \left([a,b],\eig(\Lambda_2)\right)>\delta$.
Therefore, $\tan{\theta_i} \leq \left( {l_{ii}^{(21)}}/{\delta}\right)$.

By summing over first $j$ diagonal values, we get 
$\sum_{i=1}^{j} \left( {l_{ii}^{(21)}}\right) \geq \delta \sum_{i=1}^{j} \tan{\theta_i} = \delta \cdot \Tr \left(\tan{\Theta}\right)$,
where $\Tr$ is the trace of a matrix. Using \cite{stewart1990matrix}, the sum of singular values $\sigma_1, \cdots , \sigma_j$ of a leading principal submatrix $\mathcal{L}_j$ is greater than or equal to the trace of $\mathcal{L}_j$.  Therefore, we can write for all $j$,
$\sum_{i=1}^{j} \sigma_i \geq \sum_{i=1}^{j} \left( {l_{ii}^{(21)}}\right) \geq \delta \cdot \Tr \left(\tan{\Theta}\right) = \delta \cdot ||\tan{\Theta}||_{\Phi_j},
\label{Eq_submatrix_of_R}$
where $\sigma_i$ are the singular values of the residual matrix $R$. According to Eq.~\ref{Eq_Equation_Residual}, we can use singular values of a matrix to induce matrix norms. Such norms are called the symmetric gauge functions. They are defined as follows:
\begin{align}
||R||_{\Phi_j} = \Phi_j(R) = \max_{1 \leq i_1< \cdots < i_j \leq n} \{ |\sigma_{i_1}| + \cdots +  |\sigma_{i_j}| \}
\label{Eq_norm_symmetric_gauge_function}
\end{align}

It is clear that the sum of the first $j$ singular values of matrix $R$, $\left( \sum_{i=1}^{j} \sigma_i \right)$ is less or equal to the sum of the $j$ largest singular values of $R$, $\left( ||R||_{\Phi_j} \right)$. Therefore, we have
\begin{align}
||R||_{\Phi_j} \geq \sum_{i=1}^{j} \sigma_i \geq \delta \cdot ||\tan{\Theta}||_{\Phi_j}
\label{Eq_norm_relation_SGF}
\end{align}


So far, we have Eq.~\ref{Eq_norm_relation_SGF}, which is valid for norms of the form defined by Eq.~\ref{Eq_norm_symmetric_gauge_function}. \textcolor{black}{We use Fan's theorem to generalize this result to all unitarily invariant norms, including the matrix $2$-norm and the Frobenius norm.} The Fan's theorem states that if $||R||_{\Phi_j}$ is greater than $\delta \cdot ||\tan{\Theta}||_{\Phi_j}$ for all $j$, then $||R||$ will be greater than $\delta \cdot ||\tan{\Theta}||$ for all unitarily invariant norms. Therefore, we can write 
$||R|| \geq \delta \cdot ||\tan{\Theta}||$, which \textcolor{black}{is the same} as $||\tan {\Theta} [\mathcal{R}(\tilde{X}_1),\mathcal{R}(X_1)] || \leq \frac{||R||}{\delta}$.
\end{IEEEproof}










\section{Proof of Theorem 9}
\label{Appendix_Proof_of_Theorem_9}

The gap between $\lambda$ and the set $\Lambda$ (eigengap) is larger than the gap between the eigenvalues of the perturbed space $\eig(\tilde{P}_1)$.

\begin{IEEEproof}
From Eq.~\ref{Eq_L_11} in the Appendix section, $\mathcal{L}_{11} = \lambda \Gamma^2 + \Sigma \Lambda_2 \Sigma$, we see that $\mathcal{L}_{11}$ is a diagonal submatrix \textcolor{black}{since $\Gamma$, $\Sigma$, and $\Lambda_2$ are diagonal.} The \textcolor{black}{$i$}-th diagonal entry is given as below
$l^{(11)}_{ii} = \lambda \cos^2{\theta_i} + \lambda^{(22)}_{ii} \sin^2{\theta_i}$.
After some manipulation, we obtain
$\lambda - \lambda^{(22)}_{ii} = \frac{l^{(11)}_{ii} - \lambda^{(22)}_{ii}}{\cos^2 {\theta_i}} = \frac{\delta}{\cos^2 {\theta_i}} \geq \delta $.
\end{IEEEproof}









\end{document}